\documentclass{article}
\usepackage{amssymb}
\usepackage{latexsym}
\newtheorem{thm}{Theorem}[subsection]
\newtheorem{lem}{Lemma}[subsection]
\newtheorem{cor}{Corollary}[subsection]

\newtheorem{dfn}{Definition}[subsection]
\begin{document}
\title {Calibrated Fibrations }
\author{Edward Goldstein}
\maketitle

\renewcommand{\abstractname}{Abstract}
\begin{abstract}
In this paper we investigate the geometry of Calibrated submanifolds and study
relations between their moduli-space and geometry of the ambient manifold.
In particular for a Calabi-Yau manifold we define Special Lagrangian
submanifolds
for any Kahler metric on it. We show that for a choice of Kahler metric the
Borcea-Voisin threefold has a fibration structure with generic fiber being a 
Special Lagrangian torus. Moreover we construct a mirror to this fibration.
Also for any closed $G_2$ form on a 7-manifold we study coassociative 
submanifolds and we show that one  example of a $G_2$ manifold
constructed by Joyce in \cite{Joy} is a fibration with generic fiber being
a coassociative 4-torus. Similarly we construct a mirror to this fibration.
\end{abstract}
 
\section{ Introduction } 

In their seminal paper \cite{HL} Harvey and Lawson defined the notion of 
calibrations. Let $M$ be a Riemannian manifold  and $\varphi$ be a closed
k-form. We say that $\varphi$ is a calibration if for any k-dimensional
plane $\kappa$ in the tangent bundle of $M$, we have 
\( \varphi|_{\kappa} \leq {vol({\kappa})} \).
We call a k-dimensional submanifold $L$ a calibrated submanifold if 
\( \varphi|_{L} = vol(L) \). 
It is easy to see that calibrated submanifolds minimize volume in their
homology class and thus provide examples of minimal varieties. 
For a thorough discussion and numerous examples we refer the reader to
\cite{Hv},\cite{HL},\cite{Mc}.

In this paper we study the geometry of Calibrated submanifolds and investigate
relations between their moduli-space and geometry of the ambient manifold. The
paper is organized as follows:

In section 2 we prove several comparison theorems for the volume of small balls in a Calibrated submanifold of a Riemannian manifold $M$, whose
sectional curvature is bounded from above by some $K$. Let $L$ be a minimal submanifold in $M$, $p \in L$ a point and $B(p,r)$ is a ball of radius $r$ around $p$
in $M$. Then there are a number of results on comparison between the volume of
$ L \bigcap B(p,r)$ and a volume of a ball of radius $r$ in a space form of 
constant curvature $K$ (see \cite{And}, \cite{Kim}, \cite{MY}).
Our main result in section 2 is Theorem
2.0.2, which states that if $L$ is calibrated then the volume of a ball of 
radius $r$ in the induced metric on $L$ (which is smaller than $L \bigcap 
B(p,r)$)
is greater than the volume of a
ball of the same radius in a space form of constant curvature $K$ for $r \leq r_0$ with $r_0$ depending only on the ambient manifold $M$.
As a corollary we deduce that there is an upper bound on a 
diameter of a calibrated submanifold in a given homology class.

In section 3 we investigate Special Lagrangian Geometry on a Calabi-Yau
manifold. In section 3.1 we define Special Lagrangian submanifolds for any
choice of Kahler metric on a Calabi-Yau manifold
and give basic facts pertaining to Special Lagrangian (SLag) Geometry. We will also prove a result about finite group actions on Calabi-Yau manifolds
and construct several new examples of SLag submanifolds.

In section 3.2 we study connections between moduli-space of Special Lagrangian submanifolds and global geometry of the ambient Calabi-Yau manifold. 
We will be 
interested in submanifolds, which satisfy condition $\star$ on their cohomology ring (defined in section 3.2). In particular tori satisfy this condition.
We state 2
conditions on an ambient manifold for each of those the moduli-space is not
compact. 
These conditions hold in many examples, and so we got a non-compactness theorem for the moduli-space.

In section 3.3 we use results of 2 previous sections to investigate
a Borcea-Voisin threefold in detail. We find a 
Kahler metric on it for which we can completely characterize singular 
Special Lagrangian submanifolds 
(they will be a product of a circle with a cusp curve).
Moreover SLag submanifolds don't intersect and  the compactified moduli-space 
fills the whole Calabi-Yau manifold, i.e. Borcea-Voisin 3-fold fibers with 
generic fiber being a Special Lagrangian torus. We also construct a mirror
to this fibration. 
Thus the SYZ conjecture (see \cite{SYZ}) holds in this example. 

In section 3.4 we will examine holomorphic functions on a Calabi-Yau
manifold in a neighbourhood of a Special Lagrangian submanifold. An 
immediate consequence of the fact that SLag submanifolds 
are 'Special'
is Theorem 3.4.1, which states that the integral of a
holomorphic function over SLag submanifolds is a constant
function on their moduli-space. This will give a restriction on how a
family of SLag submanifolds might approach a singularity (Corollary
3.4.1) and also will give a restriction on  SLag submanifolds asymptotic
to a cone in $\mathbb{C}^n$ (Theorem 3.4.2).

In section 4 we study coassociative submanifolds on $G_2$ manifolds. We extend
a coassociative condition for any choice of a closed (but not necessarily
co-closed) $G_2$ form. Deformation of coassociative submanifolds will still be
unobstructed and the moduli-space is smooth of dimension $b_2^+(L)$, there $L$
is a coassociative submanifold. We will show that one example of a $G_2$
manifold constructed by Joyce in \cite{Joy} is a fibration  with generic fiber
being a coassociative 4-torus. We also construct a mirror to this fibration.

There are a number of natural questions that arise from this paper. One is to
give a systematic way to construct fibrations on resolutions of torus quotients
(both for SLag and coassociative geometry). Another point is that we produced
those fibrations for certain special choices of structures on the ambient manifold (a certain choice
of Kahler metric or a certain closed $G_2$ form). We would like to get
fibrations for any other isotopic structure. If we have a 1-parameter family of
structures then we obtain a 1-parameter family of moduli-spaces $\Phi_t$. 
Suppose that $\Phi_0$ compactifies to a fibration of the ambient manifold. We
conjecture that so does each $\Phi_t$ (both for SLag and coassociative
geometries). This in particular would imply the existence of SLag fibration for
the Calabi-Yau metric on a Borcea-Voisin threefold and coassociative fibration for a parallel $G_2$ structure.
We hope to address those issues in a future paper.

{\bf Acknowledgments} : This paper is written towards author's Ph.D. at
the
Massachusetts Institute of Technology. The author wants to thank his
advisor, Tom Mrowka, for initiating him into the subject and for
constant encouragement. He is also grateful to  Gang Tian for a
number of useful conversations. Special thanks go to Grisha Mikhalkin for
explaining the Viro construction in Real Algebraic Geometry, which was
used to construct examples of real quintics.

\section {Volume Comparison for Calibrated Submanifolds } 
 If a Riemannian manifold $M$ has an upper bound $K$ on it's sectional
curvature then the volume of a sufficiently small ball in $M$ is greater
then the volume of a ball of the same radius in a space form of curvature
$K$. It turns out that this holds more generally for calibrated
submanifolds of $M$. Namely we have a following theorem :
\begin{thm}
: Let $\varphi$ be a calibrating k-form on an ambient manifold $M$ and $L$ be
a 
calibrated submanifold. Let the sectional curvature of $M$ be bounded
from above by $K$ . Let 
\(r \leq { min(injrad(M),\frac{\pi}{\sqrt{K}}) } \).              
Let $p \in L$ and $B(p,r)$ be a ball of radius r around $p$ in $M$  and \(B^K(r)\) be a
ball of radius $r$ in a k-dimensional space of constant sectional 
curvature $K$. Then 
\[ vol(L\bigcap{B(p,r)} \geq vol(B^K(r)) \]
\end{thm}
Remark : If $\varphi$ is a volume form on $M$, then this is Gunther's
volume comparison theorem.\\
{\bf Proof of Theorem 2.0.1}: The proof is based on a
following Lemma, which is
a counterpart to Rauch comparison theorem :
\begin{lem}
: Let $M$ be a (complete) Riemannian manifold whose sectional
curvature is bounded from above by $K$ and $\gamma$ : 
\( [0,t] \mapsto M \) be a unit speed geodesic. Let $Y$ be a Jacobi
field along \( \gamma  \)  which vanishes at 0, orthogonal to 
\( \gamma ' \) and \( t \leq \frac {\pi}{\sqrt{K}} \).  
Then it's length \( |Y(\theta)| \) satisfies the following differential 
inequality
\( {|Y|''} + K|Y| \geq 0\). \\
Moreover if a function \( \Psi \) is a solution to 
\( \Psi   '' + K \cdot \Psi  = 0 \) , \( \Psi(0) = 0 \)
and \( \Psi(t) = |Y(t)| \) then \[ \Psi (\theta) \geq |Y(\theta )| \]
for \( 0 \leq \theta \leq t \)  
\end{lem}
{\bf Proof }:\\
First a condition on $t$ means that $Y$ doesn't vanish on \( (0,t] \) by Rauch
Comparison theorem.
We have 
\(|Y| = \sqrt{\langle {Y,Y} \rangle } \), \(|Y|' = \frac{\langle
\nabla_{t} Y,Y  \rangle }{|Y|} \), 
\[|Y|''= \frac{| \nabla_{t} Y|^{2} - \langle Y,R(\gamma ' ,Y) \gamma '
\rangle } {|Y|} - \frac{\langle \nabla_{t} Y, Y \rangle ^{2} }{|Y|^3}   
\geq \frac{ | \nabla_{t} Y|^{2} |Y|^{2} - \langle \nabla_{t} Y ,Y
\rangle ^{2} }{|Y|^{3} } - K|Y|  
\geq - K|Y|\]
by Cauchy-Schwartz inequality. Here $R$ is a curvature operator, 
\( \gamma ' \) is a (unit length) tangent field to geodesic \( \gamma \). 
Since $Y$ is orthogonal to \( \gamma ' \) then 
\( \frac {\langle R( \gamma ' , Y) \gamma ', Y \rangle} {|Y|^{2} } \) is 
the sectional curvature of a plane through $Y$ and
\( \gamma ' \), which is less then $K$.\\
For the second claim consider \( F = \frac{|Y|}{\Psi} \). 
$\Psi$ is
positive on the interval \( (o,t] \) and hence $F$ is well defined on that
interval.
Also \( F ' = \frac{|Y| ' \Psi - \Psi ' |Y| }{\Psi ^{2} } \).\\
Consider \( G= |Y| ' \Psi - \Psi ' |Y| \). \( G(0) = 0 \) , 
\( G ' = |Y| '' \Psi - \Psi '' |Y|  \geq 0 \).\\ 
So \( G \geq 0 \) , i.e. \( F ' \geq 0 \) . Now \( F( t) = 1 \) , so 
\( F \leq 1\) i.e. \( |Y| \leq F \)            Q.E.D.\\

\noindent
Now we can prove {\bf Theorem 2.0.1}:\\
Let $d_p$ be a distance function to $p$ on $M$. Then for an open dense
set of full measure of values $t$, $t$ is a regular value of $d_p$
restricted to $L$. Let now
\[ f(t)= vol(L \bigcap B(p,t) ) \hbox{ and } g(t) = \int_{L \bigcap B(p,t)}
|\nabla_{L} d_{p}| \] 
We also consider an analogous situation on \( \overline{L} \) - a space
form of constant curvature $K$. Then \( \overline{f} = \overline{g} \)
because \(  |\nabla \overline{d_p} | = 1 \) on \( \overline{L} \).\\
For $t$ a regular value as above we have by the co-area formula :
\begin{equation} 
 f ' (t) \geq g ' (t) ~,~
 g ' (t) = vol(S_{t})
\end{equation}
and 
\begin{equation} 
\overline{f} ' = vol( \overline{S_t}),
\end{equation}
here \( S_{t} = d_{p}^{-1}(t) \bigcap L \).\\
Consider now a map $\xi$ : \( S_{t} \times [0,t] \mapsto M \),
 \( \xi ( a, \theta ) = exp( \frac{\theta}{t} exp^{-1}(a) ) \), here
\( a \in S_{t} \) , \(\theta \in [0,t] \).
Then \( vol(\xi(S_{t} \times [0,t] )) \geq f(t) \).\\  
Indeed let $\rho$ be a $k-1$ form on \( B(p,r) \) s.t. \( d \rho = \varphi
\) (such $\rho$ exists by Poincare Lemma).
Then by the calibrating condition \[ vol(\xi(S_{t} \times [0,t] )) \geq \int_{S_{t} \times [0,t]} 
\xi^{\ast} \varphi = \int_{S_{t}} \rho = \int_{B(p,t) \bigcap L} \varphi = 
f(t) \]
Also on \( \overline{L} \) we have \( \overline{f(t)} =
vol(\xi(\overline{S_{t}} \times [0,t] )) \).\\
We need to estimate \(h(t)= vol(\xi(S_{t} \times [0,t] ) \). Let $g '$ be
the product metric on \( S_{t} \times [0,t] \). Then \( h(t) = \int_{S_{t} 
\times [0,t]} Jac(d \xi) dg ' \).\\
To estimate \( Jac(d \xi) \)  at point \( (a , \theta ) \) we take an o.n.
basis \( v_{1} \ldots v_{k-1} \) to $S_t$ at $a$. Then \( d \xi (v_{i}) \)
is a value of a Jacobi field along a unit speed geodesic \( (exp(s \cdot
\frac{exp^{-1}(a)}{t}) | s \in [0,\infty))  \) at \( s= \theta \) which is
orthogonal to this geodesic, vanishes at $0$ and those length is $1$ at
\( s=t \).\\
Let $F_{t}(\theta)$ solve \( F_{t} '' + K \cdot F_{t} = 0 \), \(
F_{t}(0)=0 \), \(F_{t}(t)=1\).\\
By Lemma 2.0.1 we have \(|d \xi(v_{i})| \leq F_{t}(\theta) \), so 
\begin{equation}
Jac(d \xi) \leq (F_{t}(\theta))^{k-1} 
\end{equation}
We can consider an analogous situation on \( \overline{L} \) and in that
case we have an equality \( Jac(d \overline{\xi}) = (F_{t}(\theta))^{k-1}\).
So 
\[\overline{f(t)} = \int_{\overline{S_t} \times [0,t]}
(F_{t}(\theta))^{k-1} = vol(\overline{S_t}) \cdot \int_{[0,t]}
(F_{t}(\theta))^{k-1} d \theta = 
(by  (2)) = \overline{f} '
(t) \cdot \alpha (t) \]
(here \( \alpha (t)= \int_{[0,t]} (F_{t}(\theta))^{k-1} d \theta \)).\\
Returning now to our calibrated submanifold we deduce from (3)
and (1) that
\(f(t) \leq f ' (t) \cdot \alpha (t) \).\\
So \( \frac{f ' (t)}{f(t)} \geq \frac{\overline{f} ' (t)}{\overline{f} 
(t)} \) , i.e. \( ln(f) ' \geq (ln( \overline{f})  - \epsilon) ' \)  for
any 
\( \epsilon > 0 \). Having $\epsilon$ fixed we can choose $t_0$ small
enough s.t. \( lnf(t_{0}) \geq ln \overline{f} (t_{0}) - \epsilon \).\\
Now \(lnf(\theta ) \) is defined for a.e. $\theta$ and is an increasing 
function on $\theta$, so
\[lnf(t) \geq lnf(t_{0}) + \int_{[t_{0},t]} lnf ' \geq ln \overline{f}
(t_{0})- \epsilon  + \int_{[t_{0},t]} (ln \overline{f}) ' = ln
\overline{f} (t) - \epsilon \]
Now $\epsilon$ was arbitrary, hence \( lnf(t) \geq ln \overline{f} (t) \)
i.e. \( f(t) \geq \overline{f} (t) \)               Q.E.D .\\

\noindent  
We wish to discuss the compactification of some moduli-space of
calibrated submanifolds in a given homology class.
If we have a moduli-space $\Phi$ we can look on it as a subspace in
the space of rectifiable currents. k-dimensional currents have a mass
norm ${\bf M}$ and and a flat norm ${\cal F}$ (see \cite{Mor}, p.42) 
\[{\bf M}(L)= sup(\int_{L}\eta | \eta \hbox{ a k-form}, 
\forall p \in M :|\eta(p)| \leq 1 ) \]
\[{\cal F}(A) = inf({\bf M}(A) + {\bf M}(B) | L=A+ \partial B) \] 
Since all the
submanifolds in $\Phi$ are closed and have the same volume, then by the
Fundamental compactness theorem (theorem 5.5 in \cite{Mor}) we have that the
closure \( \overline {\Phi} \) of $\Phi$ in the flat topology is compact.\\
Also for compact subsets of $M$ there is a Gromov-Hausdorff distance
function
\(d^{GH}\), there \[ d^{GH}(K,N)= sup_{p \in K } inf_{q \in N } d(p,q)
\]
Using Theorem 2.0.1 we get 
\begin{cor}
: There is a constant \(C=C(M,{\varphi}) \) s.t. for 
\( K,N \in {\Phi} \) 
we have \(d^{GH}(K,N) \leq { C \cdot ({\cal F}(K-N))^{\frac{1}{k+1} } } \).
\end{cor}
{\bf Proof } : Suppose \( d^{GH}(K,N) = r \). Then we have a 
point \(p \in K \) s.t. \(d(p,N) =r \).\\
It is easy to construct a nonnegative function $f$ supported in a ball 
$B(p,r)$ which is equal to $1$ on a ball $B(p,r/2)$ and s.t. \( |\nabla (f)| 
\leq \frac{const}{r} \).\\
Suppose \( K-N = A+ \partial B \). Obviously \(K-N(f \varphi) \geq 
vol(B(p,r/2) \bigcap K) \)\\
\(\geq const \cdot r^{k} \) by theorem 2.0.1\\
Also \(K-N(f \varphi) = A(f \varphi) + B(df \wedge \varphi) \leq {\bf M}(A)+
\frac{const \cdot {\bf M}(B)}{r} \leq \frac{const \cdot ({\bf M}(A)+{\bf M}
(B))}{r}\).\\
So taking the infimum we get \[ const \cdot r^{k} \leq \frac{{\cal F}(K-N)}
{r} \]
which is the statement of the Corollary.         Q.E.D. 
\\
\noindent
From that we get an immediate corollary 
\begin{cor}
: If a sequence of submanifolds $L_i$ in $\Phi$ converges to
a current $L$, then it converges to the support of $L$ in Gromov-Hausdorff
topology .
\end{cor}
We now come to the main result of this section.
We wish to strengthen Theorem 2.0.1 by an analogous result for volume of
balls of radius $r$ in the induced metric on calibrated submanifolds 
(which are smaller then the balls we considered before). We have the
following
\begin{thm}
: Let $M$,$\varphi$,$p$,$L$,$K$ and \(B^{K}(r)\)  as
in Theorem 2.0.1
and let \(d_L\) be 
a distance function to $p$ on $L$ in the induced metric on $L$.
Then for \( r \leq  min(injrad(M),R(K)) \) we have :
\[vol (  x \in{L}|d_{L}(x) \leq{r})  \geq vol(B_{K}(r)) \] and \(R(K)=\pi/ 
\sqrt{K} \) for $K$ positive.
\end{thm} 
\begin{cor}:
Let $M, \varphi$ as before. Then there is an a priori bound on a diameter of calibrated submanifolds in a given homology class $\eta$.
\end{cor}
{\bf Proof of the Corollary}: Choose some $r$ satisfying conditions of theorem
2.0.2. Let $L$ be some
calibrated submanifold
in a homology class $\eta$. Let $\Gamma$ be a maximal covering of $L$ by
disjoint balls of radius $r$. Since by theorem 2.0.2 each such ball has a
a volume at least $\epsilon$ and the volume of $L$ is \( v= [ \varphi ] 
(\eta) \), then such covering exists and the number of elements in
$\Gamma$  is at most \( N = \frac{v}{\epsilon} \). Now every point in $L$
is contained in one of the balls of radius $2r$ with the same centers as 
balls in $\Gamma$.\\
So it is easy to deduce that the diameter of $L$ is at most \( 4rN \).
 Q.E.D. \\

\noindent
{\bf Proof of Theorem 2.0.2}: We wish to use the same argument as in 
the proof of
Theorem 2.0.1 for the distance function $d_L$. The problem is that $d_L$ is
not a smooth function in the $r$-neighbourhood of $p$. But we can
still smoothen it using the following technical Lemma:
\begin{lem}:
Let $L$ be a submanifold, \( p \in L \) and $d_L$
as before. 
We can pick \( \rho > 0 \) and a \( (C^{\infty}) \) function
\( 0 \leq \nu \leq 1 \) on \( [0, \infty ) \) which is $0$ on \( (0, \rho]
\), $1$ on \( [2 \rho, \infty ) \) and nondecreasing s.t.
for any positive  $\epsilon$ there is a function \( \lambda _{\epsilon}
\)on $L$ which satisfies :

1) \( \lambda_{\epsilon} \) is \( C^{\infty} \) outside of p 

2) \( d_{L} \leq \lambda_{\epsilon} \leq d_{L}(1+ \epsilon ) \)

3) \( | \nabla \lambda_{\epsilon}| \leq 1+ \nu(d_{L}) \epsilon \) 
\end{lem}
{ \bf Proof}: Pick a positive \( \rho << injrad(L) \). Choose a function
$\kappa$ on $M$ s.t. \( \kappa=1\) on \( B(p, 2 \rho ) \) and \( \kappa =
0 \) outside of \(B(p,3 \rho) \).
Choose a nonnegative radially symmetric function $\sigma$ on $\mathbb{R}^k$ 
with
support in the unit ball which integrates to 1 and let \( \sigma_{n}(x)=
n^{k} \cdot \sigma(nx) \). Then \( \sigma_{n} \) also integrates to 1.\\
Choose a nonnegative function \( \eta \leq 1 \) , \( \eta = 0 \) on \(B(p,
\frac{5 \rho}{4} ) \) and \( \eta = 1 \) outside of \( B(p, \frac{3 \rho}
{2} ) \).

Define now \( \mu^{n} : L \mapsto R \) , \[ \mu^{n}(q)= \int_{T_{q}L}
d_{L}(exp(\theta))\sigma^{n}(\theta) d \theta \]
Here \(T_{q}L \) is the tangent bundle to $q$ at $L$. Since \( \sigma^n
\) was radially symmetric function and \( T_{q}L\) has a metric, the
expression \( \sigma^{n}(\theta) \) is well defined and also integration
takes part only on a ball of radius \( \frac{1}{n} \subset T_{q}L \).\\ 
Also it is clear that \[\mu^{n}= d_{L} + o(\frac{1}{n}) \]   
The point is that for large $n$, \(\mu^n\) is a smooth function on $L$.
Indeed let us denote by \( J(a,b) \) the Jacobian of exponential map from
$a$ that hits $b$ for $a,b$ points in $L$ that are close enough. Then
\(J(a,b) \) is a smooth function on \( (a,b) \) and we can rewrite
\[ \mu^{n}(q) = \int_{L} J(q,b)^{-1} \cdot d_{L}(b) \cdot \sigma
(exp_{b}^{-1}(q)) db \]  
and it is clear from this definition that
\( \mu^n \) is a smooth function of $q$ for $n$ large enough.\\
Also one can easily prove that \( |\mu^{n}(q_{1}) - \mu^{n}(q_{2}) | \leq 
d(q_{1},q_{2}) \cdot (1 + o(\frac{1}{n})) \), hence \[|\nabla \mu^{n} |
\leq 1 + o(\frac{1}{n}) \]
Now pick \( \epsilon > 0 \). Define 
\( \lambda_{\epsilon}^{n} = (1 + \eta \epsilon)(\kappa \cdot d_{L} + (1-
\kappa) \cdot \mu^{n}) \).\\
Then \( \lambda_{\epsilon}^{n} = d_L\) on \(B(p,\frac{3 \rho}{2}) \) and
it is smooth outside of $p$. \\
One can also directly verify that we can choose a constant C s.t. for
sufficiently large $n$, the  function \( \lambda_{\epsilon} =
\lambda_{\frac{\epsilon}{C}}^{n} \)
satisfies properties 2) and 3) as desired.     Q.E.D.\\

\noindent 
Now we can prove {\bf Theorem 2.0.2}:  We will use the fact that the function \(\alpha(t) \), defined in the proof of Theorem 2.0.1, is an increasing
function of t for \( 0 \leq t \leq \frac{\pi}{\sqrt{K}} \) for $K$ positive and
for \( 0 \leq t \leq R(K) \) for $K$ negative.

Pick $\rho$ as in Lemma 2.0.2. 
Let \( \epsilon >0 \). We will follow the lines of proof of Theorem 2.0.1
for the function \( \lambda_{\epsilon} \) instead of the distance
function.
We denote by \[ f(t)= vol(\lambda_{\epsilon}^{-1}([0,t])~,~
 S_{t} = \lambda_{\epsilon}^{-1}(t) \] 
Then conditions on \( \lambda_{\epsilon} \) and the co-area formula imply
that for a regular value $t$ we have \( f'(t) \geq \frac{vol(S_{t})}
{1+ {\epsilon} {\nu(t)}} \).\\
Also we can consider  \(A_{t} = ((a, \theta)| a \in S_{t} , 0 \leq \theta
\leq d_{p}(a) ) \) (here $d_p$ is the distance to $p$ in the ambient
manifold).
We have \( \xi : A_{t} \mapsto M \), \(\xi(a, \theta) =
exp_{M}(\frac{\theta \cdot exp^{-1}(a) }{d_{p}(a)}) \).\\
As before we will have \( f(t) \leq vol(\xi(A_{t})) \) and 
\( Jac(d \xi) \leq (F_{d_{p}(a)}(\theta))^{k-1} \) (see (3), we have the same
notations as in Theorem 2.0.1). The estimate for Jacobian is true for the following reason: Let $v_1, \ldots, v_{k-1}$ be an o.n. basis to $S_t$ at $a$. Then only the normal component of $d\xi(v_i)$ to the geodesic contributes to 
$Jac(d\xi)$. The normal component can be estimated by Lemma 2.0.1.
 
So we will have \( vol(\xi(A_{t})) \leq \int_{S_t} \alpha(d_{p}(a)) da
\leq vol(S_{t}) \cdot \alpha(t) \) (here we used
the fact that $\alpha$ is an increasing function and \( d_{p}(a) \leq
d_{L}(a) \leq \lambda_{\epsilon}(a) = t \)).\\ 
Combining all this we get \[(lnf)'(t) \geq \frac{(ln \overline{f})'(t)}{1+
\epsilon \nu(t) } = [(\frac{ln \overline{f}}{1+ \epsilon \nu })' +
\epsilon \nu '/(1+\epsilon \nu)^{2} \cdot ln(\overline{f})](t)\]
Now \( \nu(t) = 0 \) for \(t \leq \rho \) and \( \nu'(t) = 0 \) for \( t
\geq 2 \rho \) and \( ln(\overline{f}) \geq -C \) for \( 2\rho \geq t \geq \rho
 \).
So \[ (lnf)' \geq (\frac{ln \overline{f}}{1+ \epsilon \nu})' - \epsilon C' \]
i.e. \( (lnf+ \epsilon C' t)' \geq (\frac{ln \overline{f}}{1+ \epsilon
\nu})' \) \\
and for $\theta$ small we have 
 \( ln(f(\theta) + \epsilon C' \theta) \geq ln(\overline{f}(\theta)) =
\frac{ln \overline{f}(\theta)}{1+ \epsilon \nu(\theta)}) \).\\
So \(ln(f + \epsilon C't) \geq ln \overline{f}/(1+ \epsilon \nu) \).\\
Here $\epsilon$ was arbitrary and we are done.           Q.E.D. 

\section{Special Lagrangian geometry on a Calabi-Yau manifold}
\subsection{ Basic properties and examples}
Let $M^{2n}$ be a Calabi-Yau manifold, $\varphi$ a holomorphic volume
form and $\omega$ a Kahler 2-form. 
If $\omega$ is a Calabi-Yau form then \(Re({\varphi}) \) 
is a calibration (see \cite {Mc})
and calibrated submanifold $L$ can be characterized by
alternative conditions : \( \omega |_{L} = 0 \) and \( Im({\varphi})|_{L}
= 0 \).
For arbitrary Kahler form $\omega$ we can define special Lagrangian (SLag) submanifolds by
those 2 conditions. The form $\varphi$ has length $f$ with respect to the metric $\omega$ (here $f$ is a positive function on $M$).
We can conformally change the metric so that the form
\( \varphi \) will have length $\sqrt{2}^n$ with respect to the new metric
$g'$. 
Then SLag submanifolds will be Calibrated by \( Re \varphi \) with respect to 
$g'$. 
\begin{lem}
Let $L^n$ be a compact connected $n$-dimensional manifold. Then the moduli-space of SLag embeddings of $L$ into $M$ is a smooth manifold of dimension $b_1(L)$. 
\end{lem}
{\bf Proof:} The proof is a slight modification of McLean's proof for a Calabi-Yau metric (see \cite{Mc}).

Let $i:L \mapsto M$ be a (smooth) SLag embedding of $L$ into $M$. Locally the moduli-space $\Gamma$ of $C^{2,\alpha}$-embeddings of $L$ into $M$ (modulo the diffeomorphisms of $L$) can be identified with the $C^{2,\alpha}$ sections of the normal bundle of $i(L)$ to $M$ via the exponential map. Also the normal bundle is naturally isomorphic to the cotangent bundle of $L$ via the map $v \mapsto i_v \omega$. Hence the tangent bundle to $\Gamma$ can be identified with $C^{2,\alpha}$ 1-forms on $L$. Let $V_k$ be the vector space of exact $C^{1,\alpha}$ $k$-forms on $L$ and let $V=V_2 \oplus V_n$. There is locally a map $\sigma: \Gamma \mapsto V$, given at an embedding $j(L) \in \Phi$ by $(j^{\ast}(\omega),j^{\ast}(Im\varphi))$. The moduli-space $\Phi$ of SLag embeddings is just the zero set of $\sigma$. Now the differential of $\sigma$ at $i(L)$ in the direction of $\alpha$ (there $\alpha$ is a $C^{2,\alpha}$ 1-form on L as above) is
\[ (d\alpha,d (f \ast \alpha)) \]
there $f$ is the length of $\varphi$ in the metric defined by $\omega$. For
$\omega$ a Calabi-Yau metric $f$ is constant. We claim that the differential is surjective and the tangent space to $\Phi$ is naturally isomorphic to the
first cohomology $H^1(L,\mathbb{R})$. To prove this consider first an operator $P$ from the space of $C^{3,\alpha}$ functions on $L$ to $C^{1,\alpha}$ $n$-forms on $L$, $P(h)=d (f\ast dh)$. We claim that $P$ is surjective onto the space of exact $n$-forms and the kernel of $P$ is the space of constant functions on $L$. Since $f$ is non-vanishing $P$ is elliptic. So to prove the surjectivity it is enough to show that the co-kernel of $P$ consists of constant multiples of the volume form on $L$. Let $\mu$ be in the co-kernel of $P$. Let $h=\ast \mu$. One easily computes that \[\int_{L}Ph \cdot \mu= \pm \int_{L} f|d^{\ast}(\mu)|^2 \]
So $d^{\ast}(\mu)=0$, hence $\mu$ is a constant multiple of the volume form on $M$. Let now $h$ be in the kernel of $P$. Then arguing as before we get that $\mu=\ast h$ is a constant multiple of the volume form on $M$, i.e. $h$ is a constant. 

Now we can prove the lemma . First we prove that $d\sigma$ is surjective. Let $\alpha$ be an exact 2-form on $L$, and $\beta$ be an exact $n$-form on $L$. We need to find a 1-form $\gamma$ on $L$ s.t. \[d\gamma=\alpha ~ , ~ d(f \ast \gamma)=\beta \]
Since $\alpha$ is exact there is a 1-form $\gamma'$ s.t. $d\gamma'=\alpha$. We are looking for $\gamma$ of the form $\gamma=\gamma'+dh$ for a function $h$. Since the operator $P$ was surjective onto the space of exact $n$-forms on $L$, we get that such $h$ exists, so $d\sigma$ is surjective, hence $\Phi$ is smooth. Next we prove that $dim(\Phi)=b_1(L)$. Let $W=ker(d\sigma)$. $W$ is the tangent space to $\Phi$ at $i(L)$. Since $W$ is represented by closed 1-forms on $L$, there is a natural map $\xi:W \mapsto H^1(L,\mathbb{R})$. We claim that this map is an isomorphism. Indeed let $a \in H^1(L,\mathbb{R})$ and let $\gamma'$ be a closed 1-form on $L$ representing the class $a$. From the properties of operator $P$ it is clear that there is a unique exact 1-form $\gamma''=dh$ s.t. $\gamma=\gamma'+\gamma''$ is in the kernel of $\sigma$. Hence $\xi$ is an isomorphism Q.E.D.

Remark: A more general setup of deformations of SLag submanifolds in a symplectic manifold was considered by S. Salur in \cite{Sem}.

In all subsequent discussions
the moduli-space will be connected (i.e. we take a connected component of the moduli-space of SLag embeddings of a given manifold $L$).

We can also define
$\Phi'$ - a moduli-space as special Lagrangian
embeddings of a given manifold $L$ into $M$ over Diff' 
(diffeomorphisms of $L$ which induce identity map on the homology of $L$).
Then \( \Phi' \) is a covering space of $\Phi$.
Now any element $\alpha$ in the first homology of $L$ induces a 1-form
\( h^{\alpha} \)on \(\Phi' \) in the following way : Let \(\xi \in \Phi '
\) and let \( L_{\xi} \) be it's support in $M$. If $v$ is a tangent
vector to \( \Phi ' \) at \( \xi \) then we can view $v$ as a closed
1-form on \( L_{\xi} \). From definition of \( \Phi ' \) it is clear
that the element \( \alpha \) induces a well defined element in \(
H_{1}(L_{\xi}) \), which we will also call \( \alpha \). So we
define \( h^{\alpha}(v) = [v](\alpha) \). Hitchin
effectively proved in \cite{Hit} that \( h^{\alpha} \) is a closed form on \(
\Phi ' \) (his notations are somewhat different from ours). Thus if we
pick \( \alpha_{1} \ldots \alpha_{k} \) a basis for the first homology of $L$ then
we have a frame of closed forms \(h^{1} \ldots h^{k} \) and
correspondingly a dual frame of commuting vector fields \( h_{1} \ldots
h_{k} \) on \( \Phi ' \). Hence any compact connected component $\Gamma$ of $\Phi$ must be a torus. Indeed the flow by commuting vector fields $h_i$ induces a
transitive $\mathbb{R}^k$ action on $\Gamma$ with stabilizer being a discrete
subgroup $A$, hence $\Gamma$ is diffeomorphic to $\mathbb{R}^k/A$ - a k-torus. 

Next we investigate finite group actions on Calabi-Yau manifolds. Suppose
that a group $G$ acts by structure preserving diffeomorphisms on $M$. We have the following
\begin{lem} 
: Suppose a SLag submanifold $L$ is invariant under the $G$ action and $G$ acts
trivially on the first cohomology of $L$. Then $G$ leaves invariant every
element in the moduli-space $\Phi$ through $L$.\\
Moreover, suppose that $g \in G$ and $x \in M-L$ in an isolated fixed point of
$g$. Then $x$ cannot be contained in any element of $\Phi$
\end{lem}
{\bf Proof} : Since $G$ is structure preserving, it sends SLag submanifolds
to SLag submanifolds. Since it leaves $L$ invariant, it preserves $\Phi$
(which is a connected component of $L$ in the moduli-space of SLag submanifolds).
From the identification of the tangent space of $\Phi$ at $L$
with $H^1(L,\mathbb{R})$ and the fact that $G$ acts trivially on
$H^1(L,\mathbb{R})$ we deduce that $G$ acts trivially on the tangent space
to $\Phi$ at $L$. Hence $G$ acts trivially on $\Phi$, i.e. it leaves each
element of $\Phi$ invariant.

To prove the second statement, consider a set $S$ of those elements in $\Phi$
which contain $x$. Obviously $S$ is closed and doesn't contain $L$.
We prove that $S$ is open and then it will be empty.

Let $L' \in S$. Any element $L''$ close to $L'$ can be viewed uniquely
as an image $exp(v)$, there $v$ is a normal vector field to $L'$. Suppose
$v(x) \neq 0$. Since $L''$ is $g$-invariant then $exp(g_{\ast}v(x))$ is
also in $L''$, there $g_{\ast}$ is a differential of $g$ at $x$. Since $L'$
is $g$-invariant then $g_{\ast}$ preserves the tangent space to $L'$ at
$x$, hence it preserves the normal space. Also since $x$ is isolated then
$g_{\ast}$ has no nonzero invariant vectors. Hence $v(x) \neq g_{\ast}v(x)$ in the normal bundle. 

Since exponential map is a diffeomorphism from a small neighbourhood of the 
normal bundle of $L'$ to $M$ we see that $exp(g_{\ast}v(x))$ is not in $L''$-
a contradiction. So $v(x)=0$ i.e. $L'' \in S$      Q.E.D.

As for examples of special Lagrangian submanifolds, many come
from the following setup : Let $M$ be a Calabi-Yau manifold and $\sigma$
an antiholomorphic involution. Suppose $\sigma$ reverses $\omega$.
Then the fixed-point set of $\sigma$ is a special Lagrangian submanifold.
For a Calabi-Yau metric $\omega$ the condition $\sigma$ reverses $\omega$ is
equivalent to $\sigma$ reversing the cohomology class $[\omega]$, which often 
can be easily verified.
Indeed suppose $\sigma$ reverses $[\omega]$. Then \( -\sigma^{\ast}(\omega)\)
is easily seen to define a Kahler form, which lies in the same cohomology class
as $\omega$ and the metric it induces is obviously equal to \(\sigma^{\ast}(g)
\), i.e. it is Ricci-flat. Hence by Yau's fundamental result (see \cite{Yau}) 
we have \(-\sigma^{\ast}(\omega)= \omega \).

We wish to discuss 2 collections of such examples. In both cases $M$ is a
projective manifold defined as a zero
set of a collection of real polynomials.
Then the conjugation of the projective space induces an anti-holomorphic
involution which reverses the Fubini-Study Kahler form, hence it also reverses the Calabi-Yau form in the same cohomology class.
The fixed point set is a submanifold of a real projective space .\\

\noindent
 Our first example with be a complete intersection of hypersurfaces of
degree 4 and 2 in \( \mathbb {C}P^5 \).\\
First we note that a 2-torus can be represented as surface of degree 4 in
$\mathbb{R}^3$. Indeed a torus can be viewed as a circle bundle over a circle
\( ((x,y,0)|x^{2} + y^{2} = 1 ) \) in $\mathbb{R}^3$, there a fiber over a 
point
\( a=(x,y,0) \) is a circle of radius \( \frac{1}{2} \) centered at $a$
and passing in a plane through \( a , (0,0,1) \) and the origin . 
If \( (x,y,z) \) is a point on our torus then it's distance to a point
\( (\frac{x}{\sqrt{x^{2} + y^{2}}}, \frac{y}{\sqrt{x^{2}+y^{2}}},0) \) is
\( \frac{1}{2} \).\\
If we compute we get \( 1+x^{2} + y^{2}+z^{2} -2 \sqrt{x^{2}+y^{2}} =
\frac{1}{4} \) , i.e. 
\( p(x,y,z) = ( \frac{3}{4} + x^{2} + y^{2}+z^{2})^{2} - 4(x^{2}+y^{2}) =
0 \).

So in inhomogeneous coordinates \( x_{1} \ldots x_{5} \) on \(\mathbb{R}P^5\) 
the zero locus of 2 polynomials \(p(x_{1},x_{2},x_{3}) \) and
\(q(x_{4},x_{5}) \) (there \( q(x,y)=x^{2}+y^{2}-1 \) ) is a 3-torus.\\
If we consider the corresponding homogeneous polynomials on \(\mathbb{R}P^5\),
then it is easy to see that there is no solution for \(x_{6}=0 \).
So the zero locus of those polynomials in \(\mathbb{R}P^5\) is a 3-torus.
If we perturb them slightly so that the corresponding complex 3-fold in
\(\mathbb{C}P^5 \) will be smooth then we obtain the desired example.\\

\noindent
 Other examples are quintics with real coefficients  in \(\mathbb{C}P^4\). In 
that case real quintics would be  special Lagrangian submanifolds. R.Bryant
constructed in \cite{Br} a real quintic, which is a 3-torus .\\
We will construct, using Viro`s technique in real algebraic geometry 
(see \cite{Vir}), real
quintics $L_k$ which are diffeomorphic to projective space
\(\mathbb{R}P^3\) with k 1-handles attached for \( k = 0, 1 , 2 , 3\).  If
$k=3$ then \(b_{1}(L_{3})=3 \) and the cup product in the first cohomology
of $L_3$ is $0$.

The construction goes as follows: First in inhomogeneous coordinates
\( x_{1} , \ldots ,x_{4} \) on \(\mathbb{R}P^4\) we consider a polynomial 
\(p \cdot
q \), there \( p(x_{1}, \ldots , x_{4} ) = x_{1}^{2} + \ldots + x_{4}^{2}
- 1 \) and \( q = x_{4} \). In \(\mathbb{R}P^4\) the zero locus of the 
polynomial will be
\( \mathbb{R}P^3 \bigcup S^3 \), there $\mathbb{R}P^3$ is a zero locus of $q$,
$S^3$ is a zero locus of $p$ 
and $\mathbb{R}P^3$ intersects $S^3$ along 
a 2-sphere \( S^{2} \subset \mathbb{R}^3 = (x_4=0)\).
 
Now we consider \( f = pq -\epsilon h \), there $h$ is some polynomial
of degree up to 5 and $\epsilon>0$ is small enough.
The Hessian of $pq$ on $S^2$
is nondegenerate along the normal bundle to $S^2$ and vanishes along the axes 
of the normal bundle (the axes are a normal bundle of $S^3$ to $S^2$ and of
$\mathbb{R}P^3$ to $S^2$). If we look on those axes locally as coordinate
axes then $pq= xy$ in those coordinates.

Suppose first $h$ is non-zero along $S^2$. We can assume that $h>0$ on $S^2$.
The zero locus of $f=pq - \epsilon h$ will live in 2 quadrants in which
$xy>0$. Thus in zero locus of $f$, 
a part $A_1$ of $\mathbb{R}P^3$ outside $S^2$ will "connect" with one 
hemisphere $S'$ of
$S^3$, and the part $A_2$ inside $S^2$ will connect with
the other hemisphere $S''$ and the zero locus of $f$ is a disjoint union of 
$\mathbb{R}P^3$ and $S^3$.\\
Suppose the zero set of $h$ intersects our $S^2$ transversally along  k 
circles such that no circle will lie in the interior of the other. We can 
assume w.l.o.g. that on the exterior $V$ of those circles $h$ is positive.
Then along $V$, $A_1$ connects to $S'$ and $A_2$ connects to $S''$ as before.
Along the interior of every circle, $A_1$ connects to $S''$ and $A_2$ to
$S'$. So near interior of these circles we get 1-handles connecting 
$\mathbb{R}P^3$ with $S^3$.
So the zero locus of $f$ will be $\mathbb{R}P^3$ 
and $S^3$
connected by k 1-handles, i.e. it will be an $\mathbb{R}P^3$ with $k-1$ 
1-handles attached.\\
It is not hard to find examples of such $h$ for small values of k.
For instance for $k=4$ (i.e. for $L_3$) we can take \[h= ((x_1-1/3)^2 +
(x_2 -1/3)^2 - 1/16)((x_1+1/3)^2 + (x_2 +1/3)^2 - 1/16) \] and the zero
locus of $h$ intersects $S^2 \subset \mathbb{R}^3$ in 4 circles.

\subsection{ Non-compactness of the moduli-space}
In this section we will consider connections between the moduli-space of
SLag submanifolds and global geometry of the ambient Calabi-Yau
manifold $M$.\\
Let $\Phi$ be the moduli-space of SLag submanifolds. We have a fiber
bundle $F$ over $\Phi$, \( F \subset M \times \Phi \),
\( F =((a, L)|a
\in M , L \in \Phi \) s.t. $a \in L$) .\\
We have a natural projection map \( pr : F \mapsto \Phi \),
whose fiber is the support of the element in $\Phi$ and the evaluation map
\( ev : F \mapsto M \), \( ev(a, L)=a \).\\
Also the tangent space to a point \( (a , L) \in F \) naturally splits
as
\( T_{a}L \oplus T' \), there \(T_{a}L\) is
the tangent space to $L$ at $a$ (a tangent space to the fiber) and 
\(T'= ((v(a),v)|v \) is a variation v. field to the moduli-space
and \(a(v) \) is the value of $v$ at $a$).

Let $L$ be a compact k-dimensional oriented manifold with \(b_{1}(L)=
k\). We say that $L$ satisfies condition $\star$ if for $\alpha_1$ \ldots
$\alpha_k$  a basis for \( H^{1}(L) \) we have 
\( \alpha_{1} \cup \) \ldots \( \cup \alpha_k \neq 0 \) . 
This holds e.g. if $L$ is a torus. On the other hand the real quintic
with \( b_{1} =3 \) that we constructed in section 3.1 doesn't satisfy
condition $\star$.          
\begin{thm}
: Let $L$ be a special Lagrangian submanifold with \(b_{1}(L)
= k\).\\
Suppose $L$ satisfies $\star$ .
Suppose some connected component of \( \Phi ' \) is compact. Then the 
Betti numbers of $M$ satisfy : \(b_{i}(M) \leq b_{i}(L \times T^{k}) \)
(here $T^k$ is a $k$-torus).\\
Suppose we have $G,g,x$ satisfying conditions of Lemma 3.1.1. Then $\Phi$ 
itself is not compact.
\end{thm}
{\bf Proof of Theorem 3.2.1}: Suppose $L$ satisfies $\star$. First prove that 
$\Phi$ is orientable (in fact it has a natural volume element $\sigma$). 
Let \( L' \in \Phi \) and \( v_{1} \ldots v_{k} \) be 
elements of the tangent space to $\Phi$ at $L'$. So \( v_{1} \ldots
v_{k} \) are closed 1-forms on $L'$ and we define:
\[\sigma( v_{1} \ldots v_{k} ) = [v_{1}] \cup \ldots \cup [v_{k}]
(L')\] 
Suppose that $\Phi$ is compact, or a connected component $\Gamma$ of $\Phi'$
is compact. In each case we have an evaluation map as before. We will prove that in both cases it has a positive degree. We will give a proof for $\Phi$, the
proof for $\Gamma$ is analogous.

First $\Phi$ has a natural volume element $\sigma$ described above.
So the $2k$-form \( \alpha = pr^{\ast}(\sigma) \wedge ev^{\ast}(Re
\varphi) \) is the volume form on $F$.\\ 
Let \( L_{\phi} \in \Phi \) and  \( \alpha_1
\ldots \alpha_k \) be a basis for \( H^{1}(L_{\phi}) \) s.t. \(
\alpha_1 \cup \ldots \cup \alpha_k [L_{\phi}] = 1 \). Then we can consider
corresponding
vector fields \(v_{1} \ldots v_{k} \) along $L_{\phi}$, which
form a frame for the bundle $T'$ (described in the beginning of this
section) restricted to \( L_{\phi} \). So $[i_{v_j}\omega]= \alpha_j$ and
\(pr^{\ast}(\sigma)(v_{1}, \ldots , v_{k} ) =1 \).

Let now $\eta$ be a Riemannian volume form on $M$. Then we have \[deg(ev) = 
\int_{F}ev^{\ast}(\eta) / vol(M) \] Since $F$ is a fiber bundle we
can use integration over the fiber formula to compute:
\[ \int_{F}ev^{\ast}(\eta) = \int_{\Phi}(\int_{L_{\phi}} i_{v_1} \ldots i_{v_k}
ev^{\ast}(\eta) )d \phi \] ( of course we choose \( \alpha_1 \ldots \alpha_k
\) for each \( L_{\phi} \) ). \\
Also \( i_{v_1} \ldots i_{v_k}
ev^{\ast}(\eta) \) is  easily seen to be equal to \(i_{v_1} \omega
\wedge \ldots \wedge i_{v_k} \omega \) (all restricted to the fiber
\(L_{\xi} \) ).\\
So  \( \int_{L_{\xi}} i_{v_1} \ldots i_{v_k} ev^{\ast}(\eta) = \alpha_1
\cup \ldots \cup \alpha_k (L_{\phi}) = 1 \).\\
So \( deg(ev) = \int_{\Phi} 1 /vol(M) = vol(\Phi)/vol(M) > 0 \).

Now let $\Gamma$ be compact. Let $F'$ be a corresponding fiber bundle over $\Gamma$.
First we claim that \(b_{i}(F') \geq
b_{i}(M) \).
Suppose this is not true for some $0<i<k$ . Then we have \( \beta \in
H^{i}(M) \) s.t. \(ev^{\ast}(\beta) = 0 \in H^{i}(F') \).
By Poincare duality we can find \( \alpha \in H^{k-i}(M) \) s.t. \( \alpha
\cup \beta \) is the generator of \(H^{k}(M) \). Since \(deg(ev) \neq 0
\) we have \( 0 \neq ev^{\ast}( \alpha \cup \beta ) = ev^{\ast}(\alpha)
\cup ev^{\ast}(\beta) = 0 \) - a contradiction.

Next we prove that \(b_{i}(F') \leq b_{i}(L \times T^{k}) \).
We know from Section 3.1 that $\Gamma$ is a $k$-torus. 
Also if we fix a basis for each cohomology \( H^{i}(L,\mathbb{R}) \) then
each
fiber of the fibration $F'$ over $\Gamma$ has a canonical basis for
cohomology (because $\Gamma$ is a space of SLag embeddings of a
manifold $L$ modulo diffeomorphisms which induce identity map in 
homology).

Let \( U_{\alpha} \) be a good cover of $\Gamma$. Then
\(pr^{-1}(U_{\alpha}) \) is a cover of $F'$. Consider a corresponding
double
complex (see \cite{BT}, p.96) \(K=K^{p,q}=C^{p}(pr^{-1}(U_{\alpha})_{q}) \)
(here  \(pr^{-1}(U_{\alpha})_{q} \) is a collections of all $q$
intersections of elements in \(pr^{-1}(U_{\alpha}) \)).

Then this complex computes the DeRham cohomology of $F'$. Also the
cohomology of this complex is computed by the corresponding Leray spectral
sequence (see \cite{BT}, p.165) \(E_{r}^{m} = \oplus E_{r}^{p,q} \) and the
second term is given by \(E_{2}^{p,q} = H_{\delta}^{p}(H_{d}^{q}(K)) \),
there \( H_{d}^{q} \) is a sheaf over $\Gamma$ given by \(H_{d}^{q}(U)=
H^{q}(pr^{-1}(U)) \). Since each fiber has a canonical basis for cohomology we see that \( H_{d}^{q} \) is a constant sheaf \(H^{q}(L) \) over
$\Gamma$ and hence \(E_{2}^{p,q} = H^{q}(L) \otimes H^{p}(\Gamma) \), so
\( dim(E_{2}^{m}) = b_{m}(Y \times T^{k} ) \).
So we  have \(b_i(M) \leq b_i(F')=dim(E_{\infty}^{m}) \leq dim (E_{2}^{m})=b_i(L \times T^k) \).

Finally suppose $G,g,x$ satisfy conditions of Lemma 3.1.1. Suppose $\Phi$
is compact. Then the degree of the evaluation map is positive, hence it it 
surjective. But $x$ is not in the image of the evaluation map- a contradiction.
Q.E.D.

\subsection{ Special Lagrangian fibration on a Borcea-Voisin threefold}
In this section we will use results of 2 previous sections to investigate
one example of a Calabi-Yau manifold in detail. We study the 
Borcea-Voisin threefold $M$. It will turn out that for a suitable choice
of a Kahler metric on $M$ we can prove that $M$ fibers with a generic
fiber being a Special Lagrangian torus. Moreover we construct a mirror
fibration to our Special Lagrangian fibration.

The property, which we will utilize in studying $M$, is the fact that $M$ is a 
union of neighbourhoods, each of those is biholomorphic to a product of an
elliptic curve with an open neighbourhood on a $K^3$ surface.
 This will enable us to define a certain Kahler metric on $M$ such that in any such neighbourhood SLag submanifolds in our moduli-space will look like $S^1 \times T$, there $S^1$ is a circle and $T$ is
a pseudoholomorphic 2-torus. This will enable us to use Gromov's compactness theorem to study the compactification $\overline{\Phi}$ of the moduli-space $\Phi$ of SLag tori. A crucial point will be the fact that the boundary
$\overline{\Phi} - \Phi$ has dimension 1 (i.e. co-dimension 2). So the total space $F$ (see section 3.2) compactifies to a pseudo-cycle (see \cite{RT}). Moreover both $F$ and $M$
are orientable. So we can define the degree of the evaluation map $deg(ev)$ as follows:

Let $x \in M$ be a regular point (i.e the evaluation map is transversal to $x$
and $x$ is not contained in the total space of the boundary of the moduli-space). Then $deg(ev)$ is a sum over the preimages  $ev^{-1}(x)$  of their signs.  By usual transversality this is well-defined and in our example we will have $deg(ev)=1$, i.e. the moduli-space fills the whole manifold $M$.

We will adopt the definition of $M$ from \cite{Lu} , so we define $M$ to be
the resolution of singularities of a 6-torus \(T^6= T^{2} \times T^{2}
\times T^{2} \) by \( \mathbb{Z}_2 \oplus \mathbb{Z}_2 \), there the 
generators of the
$\mathbb{Z}_2$ actions are :
\[ \alpha:z_{1} \rightarrow -z_{1} + \frac{1}{2} , z_{2} \rightarrow
-z_{2} + \frac{1}{2}, z_{3} \rightarrow z_{3} \]
and \[ \beta: z_{1} \rightarrow -z_{1},z_{2} \rightarrow z_{2} , z_{3}
\rightarrow -z_{3} \]
The fixed locus of $\alpha$ is 16 2-tori :\\
\( A \times A \times T^{2} \) \\
and the one of $\beta$ is 16 2-tori:\\
\( B \times T^{2} \times B \).\\
Here \(A \subset T^2 \) is a set \( \{ \frac {2+2i \pm 1 \pm i}{4} \}\) and
\(B \subset T^2 \) is a set \( \{ 0, \frac {1}{2}, \frac{i}{2} ,
\frac{i+1}{2} \} \).\\
The fixed loci of $\alpha$ and $\beta$ do not intersect. Also \( \alpha
\circ \beta \) has no fixed points.

Consider a fixed torus $T^2$ (say of $\alpha$). Near $T^2$ the quotient looks
like \(V= (U/ \pm 1) \times T^2 \), there $U$ is a ball of radius $r$ around
a fixed point in $T^4$ (we can also view it as a neighbourhood of the
origin in $\mathbb{C}^2$).\\
The resolution of singularity \( U/ \pm 1 \) is an $r$-neighbourhood
\( \overline{U} \) of the zero set in the total space of \(\gamma^{\otimes
2}\), there $\gamma$ is the universal line bundle over
$\mathbb{C}P^1$ (thus the singular point will be replaced by \(\mathbb{C}P^1
)\).
The resolution has a 1-parameter family of
HyperKahler metrics $\omega_t$ (the Eguchi-Hanson metrics) (see \cite{Joy}, p.304). Their Kahler potentials $f_t$ are given by \[f_t(u)= \sqrt{u^2 + t^2} + t^2 logu - t^2 log(\sqrt{u^2+t^2}+t^2) \] here $u= |z_1|^2 + |z_2|^2$.
So we replace \(V=(U/\pm 1)
\times T^{2} \) by \( \overline{V}= \overline{U} \times T^{2} \).
Now the Eguchi-Hanson Kahler form $\omega_t$ can be glued to a Euclidean metric
inside of $U/ \pm 1$ by gluing their Kahler potentials for $t$ small enough
and we obtain a Kahler
form $\omega '$ on \( \overline{U}\). We consider the
corresponding product metric in $\overline{V}$. Doing this for every fixed
2-torus (both of $\alpha$ and of $\beta$) we get a Kahler metric on
$M$, which is Euclidean outside
neighbourhoods described above.

If we take a family of 3-tori \(T_{a,b,c} \subset M \),
\[T_{a,b,c}=((z_{1},z_{2},z_{3})|Rez_{1}=a, Rez_{2}=b,Rez_{3}=c ) \]
which
don't intersect neighbourhoods of fixed components, then they will be 
SLag tori in $M$ (according to the definition in the beginning of section
3.1). 
We would like to see what happens to this family then it's
elements intersect some neighbourhood \( \overline{V}\) described above.

We wish to point out that P. Lu considered SLag submanifolds for a Calabi-Yau metric on $M$ in \cite{Lu}. He was able to produce a big open set of those submanifolds. We are using a different metric and this will allow us to characterize
the compactified moduli-space and to prove that $M$ fibers over it.

We return to the question of characterizing those elements which intersect a 
'bad' neighbourhood. If this is a neighbourhood of a fixed component of 
$\alpha$, we will consider the following setup :\\
Consider a $\mathbb{Z}_2$ action on $T^4$ with a generator 
\[ \alpha':z_{1} \rightarrow -z_{1} - \frac{1}{2},z_{2} \rightarrow -z_{2}
- \frac{1}{2} \] 
Then this action has 16 fixed points and the resolution of singularities
gives a $K^3$ surface.
Also our original manifold can be viewed as 
\[ (T^{4}/ \alpha ' \times T^{2})/ \beta \]
For each fixed point of $\alpha '$ we introduce a neighbourhood $U$ as
before. We will have to consider a bigger neighbourhood \(X= ((z_{1} + ia,
z_{2} +ib)| (z_{1},z_{2}) \in U, a,b \in \mathbb{R}/\mathbb{Z}) \) in $T^4$ and a
corresponding
neighbourhood in the quotient, which we call $X'$.
We will consider a resolution of singularities \(\overline{X}
\) (now we have 4 singular points in $X'$) and a corresponding domain
\( \overline{W} = \overline{X} \times T^{2} \) in $M$ (since \(
\beta(\overline{W}) \bigcap \overline{W} = \emptyset \),  we can view
\(\overline{W}\) as a domain in $M$).

Consider a canonical \( (2,0) \)
form \(\eta = Re \eta + iIm \eta \) on
\(X \subset T^4 \) (The collection \(w, Re \eta,Im \eta \) is the
standard HyperKahler package on $T^4$).
Then \(\eta\) lifts to a holomorphic (2,0) form on \( \overline{X} \) and 
we have \(i \cdot \eta \wedge dz_{3} = \varphi \) - a holomorphic $(3,0)$ form
on $M$.\\
Since the Calabi-Yau structure on \(\overline{W}\) is a product structure,
then in \( \overline{W}\)
SLag submanifolds, which come from a connected component of a family
\(T_{a,b,c}\), will look like :          \( L \times T_c \) \\
there \(T_{c} = (z|Rez =c) \subset T^2 \) and $L$ is a SLag submanifold 
of \( \overline{X}\) w.r. to $\omega', \eta$.
Indeed submanifolds described above form a 3-dimensional
moduli-space, which is contained in the original moduli-space, hence it
coincides with the connected component of the original moduli-space.

The package \( \omega' ,Re \eta , Im \eta \) is not a HyperKahler
package on \( \overline{X}\). We can however normalize $\omega'$ to some
\(\omega'' \)(by multiplying it by a positive function) s.t. in the metric
defined by \( \omega '' \) the form $\eta$ will have length $\sqrt{2}$. In this
case \( \omega '', Re \eta,Im \eta \) will be a (non-integrable)
HyperKahler package on \( \overline{X}\).

A SLag submanifold $L \subset \overline{X}$ is defined by conditions \(\omega' |_{L} = 0,
Im \eta |_{L} =0\), which is of course equivalent to \( \omega '' |_{L} =
0, Im \eta |_{L} = 0 \), hence $L$ is a SLag submanifold with respect to
our HyperKahler package, which is equivalent to being pseudoholomorphic
submanifold with respect to (non-integrable) almost complex structure
defined by \( Re \eta \).

Now elements in our moduli-space on $M$ have trivial self-intersection.
Since those elements, which live in \( \overline{W} \) look like \( L
\times T_a \) and $L$ is pseudoholomorphic in \( \overline{X} \), it
is clear that they do not intersect. Also to characterize the
moduli-space,
it is obviously enough to study pseudoholomorphic tori in
\(\overline{X}\). Indeed the boundary of $\overline{X}$ is fibered by pseudoholomorphic tori. Hence one easily deduces 
that an element in our moduli-space is 
either contained in $\overline{W}$ or in it's complement. 

We return now to our pseudoholomorphic tori living in $\overline{X}$.
We can view those tori as living in $K^3$ as before. All of those tori
will carry the same homology class $h$. We would like to study the boundary of our moduli-space, i.e. to understand what are possible limits of a sequence of such tori.
A limit of any sequence $T_i$ of such tori, by Gromov's compactness
theorem (see \cite{RT} or \cite{Wol}), is a cusp curve with at most 1
component being a torus
and the rest are pseudoholomorphic spheres.

Suppose we have 1 component being a torus $T$. We want to prove that $T$ is the only component, it is embedded and lives in the original moduli-space and $T_i$ converge to $T$ in the moduli-space.

We can represent $T$ as a 
composition $\alpha \circ \rho: T^{2} \mapsto K^3$, there $\rho: T^{2} 
\mapsto T^2$ is a k-fold covering and $\alpha: T^{2} \mapsto K^3$ is a simple
curve (see \cite{MS}, p.18).\\
Let $T'= \alpha(T^{2})$. Since $T$ doesn't
intersect any of the $T_i$ we have that \( [T'] \cdot h = 0 \).
Also \(h = k[T'] + \Sigma[S_{i}]\) for some pseudoholomorphic spheres $S_i$.\\
We have \( [T'] \cdot \Sigma [S_{i}] \geq 0 \) with equality iff there are
no $S_i$ (because the limiting curve is connected). Also since $\alpha$ is 
simple we get by theorem 7.3 in \cite{MW} that
$[T'] \cdot [T'] \geq 0$ with equality iff $T'$ is embedded. From all this we 
deduce that there are no rational components and $T'$ is embedded.

Since $T'$ is pseudoholomorphic, it is a SLag torus, hence it admits a 
2-dimensional deformation family of SLag tori. Moreover this family fills some neighbourhood of a point in $T'$.
Indeed let $\alpha_{1},\alpha_{2}$ be
2 generators of $H^{1}(T')$ and $v_{1},v_{2}$ are corresponding deformation
v.fields. Then $i_{v_1}\omega' \wedge i_{v_2} \omega'$ is nonzero in 
$H^{2}(T')$, so it doesn't vanish at some point $p \in T'$. 
Hence $v_1$ and 
$v_2$ are linearly independent at $p$, so the moduli space fills a 
neighbourhood $U$ of $p$. 

Also $T_i$ converge to $T$ in Gromov-Hausdorff topology, hence they intersect
$U$ for $i$ large enough. So they intersect elements in the moduli-space 
through $T'$. Now these elements are pseudoholomorphic and carry a homology 
class $h/k$ and $h \cdot h = 0$. So we deduce that $T_i$ are in the 
moduli-space through $T'$ and so $k=1$ and $T=T'$ is in the original 
moduli-space. So the boundary of the moduli-space consists of unions of 
spheres.

Let $\Sigma$ be a cusp curve on the boundary on the moduli-space. Then
as before $\Sigma$ doesn't intersect any of the smooth pseudoholomorphic tori.
We want to prove that $\Sigma$ doesn't intersect any over cusp-curve on the
boundary of the moduli-space. To prove that we will have to understand the
second homology  $H_2(\overline{X},\mathbb{Z})$ in detail. 
Suppose $X$ (see p.14) is a neighbourhood of the torus $T=(Re(z_1)=1/4, Re(z_2)=1/4)$. 
Then
in the quotient $X'$, $T$ becomes a sphere $S$ and one sees that $S$ is a 
strong deformation retract of $X'$. The resolution of singularities 
$\overline{X}$ is obtained from $X'$ by replacing the singular points with the exceptional spheres
$S_i$. Also there is a sphere $S'$ in $\overline{X}$, which projects onto
$S$. 
One can easily deduce from Mayer-Vietoris sequence that any element $\alpha$ in
$H_2(\overline{X}, \mathbb{Z})$ can be represented as $\alpha = \lambda \cdot [S']+
\Sigma\lambda_i \cdot [S_i] +$ torsion, there $\lambda, \lambda_i$ are integers. Now $\int_{S'} Re(\eta)= 1/2 \cdot \int_{T}
Re(\eta) = 1/2 $ and $\int_{S_i} Re(\eta)= 0 $. So $\int_{\alpha}Re \eta= \lambda \cdot 1/2= \lambda/2$. So if $\int_{\alpha}Re\eta > 0$ then $\int_{\alpha}Re \eta \geq 1/2$.

Let now $\Sigma$ be as before. Then $\Sigma$ represents a homology class $h$
and $\int_{h}Re(\eta)=1$. Also the integral of $Re(\eta)$ on every component 
of $\Sigma$ is at least $1/2$, so $\Sigma$ has at most 2 components.
Let $\Sigma '$ be another cusp curve on the boundary of the moduli-space.
Suppose $\Sigma '$ intersects $\Sigma$. Since $h \cdot h = 0$ we see that 
$\Sigma$ and $\Sigma '$ must have a common component. Suppose $\Sigma$ has a 
component $P$ which is not in $\Sigma'$. Then $0= [P] \cdot h = [P] \cdot
[ \Sigma '] > 0 $- a contradiction. So $\Sigma$ and $\Sigma '$ have  same
components, and since their total number (counted with multiplicity) is at 
most 2, then $\Sigma = \Sigma '$.

Finally we prove that the number of exceptional spheres is finite.
As we have seen, there are 2 types of exceptional curves:

1) A curve with 2 components $A_i$ and $B_i$. Then $0= [A_i] \cdot h=
[A_i] \cdot ([A_i] + [B_i])$. Now $[A_i] \cdot [B_i] > 0$, so $[A_i] \cdot
[A_i] < 0 $.

If $A_j,B_j$ is another curve like that, then we have seen that $A_i$
doesn't intersect it, so in particular $[A_i] \cdot [A_j] = 0$.
So one easily sees that the numbers of such curves is at most
$5= b_2(\overline{X})$. 

2) A curve with 1 component (possibly with multiplicities). Let this
curve be $k \cdot P_i$, there $P_i$ is a primitive rational curve and
$k \cdot [P_i]= h $. To 
study those $P_i$ we make following observations: There is a $\mathbb{Z}_2
\oplus \mathbb{Z}_2$ action on $T^4$ with generators

$\gamma_1: (z_1,z_2) \mapsto (z_1 + i/2, z_2)$

$\gamma_2: (z_1,z_2) \mapsto (z_1, z_2 + i/2)$ 

This action commutes with $\alpha'$ action, and hence induces an action on
$K^3$. It also preserves $\overline{X}$.

Next we find elements in $K^3$, which do not have a full orbit under the
action. A point $(z_1,z_2)$ doesn't have a full orbit if it is preserved
under one of $\gamma_1, \gamma_2, \gamma_1 \circ \gamma_2$.
Now the fixed points are:

$Fix(\gamma_1)= ((z_1,z_2): (z_1 +i/2, z_2)= (-z_1 +1/2, -z_2+1/2))$. 
These are 2 points, disjoint from exceptional spheres.    A similar
analysis for $\gamma_2$ and $\gamma_1 \circ \gamma_2$ produces 2 points
for each.

Now the actions of $\gamma_i$ are structure preserving on $\overline{X}$,
so they send SLag tori to Slag tori. Moreover they preserve an open set
of tori $Re(z_i)=const$ in our moduli-space. So by Lemma 3.1.1 they leave elements
of the moduli-space invariant. Hence they preserve the
limiting curve $P_i$ (because the convergence is in particular a
Gromov-Hausdorff convergence).

For a limiting curve $P_i$,
consider $\chi(P_i) = [P_i] \cdot [P_i] - c_1(K^3)([P_i]) + 2 = 2$.
Then by theorem 7.3 of \cite{MW} we can count $\chi(P_i)$ by
adding contributions of singular points (which are double points or branch 
points), and each singular point gives a positive contribution. So $P_i$ has
singular points and there are at most 2 of those.

Let $x$ is a singular point.
Then it's orbit under $\mathbb{Z}_2 \oplus \mathbb{Z}_2$ action consists of 
singular points. So it cannot have length 4, so $x$ is one of 6 points $D$ 
with orbit of length 2. So $P_i$ contains at least 2 points of the set $D$.

If $P_j$ is another curve of type 2, then $[P_i] \cdot [P_j] = 0$, so they
don't intersect. Also $P_j$ contains at least 2 points from the set $D$.
So it is clear that the number of $P_i$ is at most 3.
 
From all that we deduce that our moduli-space can be compactified to a pseudo-cycle.  Also any point $x$ outside of bad neighbourhoods has a unique preimage
in the smooth part of the moduli-space, so we deduce that the degree of the
evaluation map is 1, so the compactified moduli-space fills the whole manifold
$M$. Also elements of the compactified moduli-space don't intersect, so
$M$ fibers with generic fiber being a SLag torus. Also the fibration is 
smooth over the smooth part of the moduli-space. To prove that we need to prove that the differential of the evaluation map is an isomorphism. This is clearly
true outside our `bad' neighbourhoods. Inside a bad neighbourhood, it is
enough to prove that variational vector fields to our pseudoholomorphic tori
do not vanish. But this follows from a standard argument that each zero of such
a vector field gives a positive contribution to the first Chern class of the 
normal bundle, which is trivial.

We want to point out that this example, then we can use
HyperKahler trick and local product structure to study limiting SLag
submanifolds, is quite ad hoc and some new ideas are needed to study
singular SLag submanifolds in general.

Next we wish to construct a mirror, i.e. to compactify the dual fibration.
Let $  \stackrel{M}{\stackrel{\downarrow }{\overline{\Phi}}}$ be a fibration
over the compactified moduli-space and let $ \stackrel{M_0}{\stackrel{\downarrow }{\Phi}}$ be a restriction of this fibration over the (smooth)
moduli-space, 
there $M_0 \subset M$ is an open subset. Let $a \in \Phi$ and $L_a$ be a fiber.
We have a vector space $V_a=H_1(L_a,\mathbb{R})$ and a lattice $\Lambda_a=
H_1(L_a, \mathbb{Z})$ in it, and so we get a torus bundle $V_a/\Lambda_a$ 
over $\Phi$. By dualizing each $V_a$ we get a dual bundle $V_a^{\ast}/
\Lambda_a^{\ast}$. We will adopt the following definition of a topological mirror from M. Gross's paper (see \cite{MG1})

\begin{dfn}
:Let $\stackrel{M'}{\stackrel{\downarrow }{\overline{\Phi}}}$ be another fibration with $M'$ smooth and a corresponding fibration
$\stackrel{M_0'}{\stackrel{\downarrow }{\Phi}}$ is a smooth torus 
fibration. Let $V_a'/\Lambda_a'$ as before.
We say that $M'$ is a topological mirror to $M$ if there is a fiberwise linear 
isomorphism $\rho : V_a^{\ast}/ \Lambda_a^{\ast} \mapsto V_a'/ \Lambda_a'$
over $\Phi$.
\end{dfn}
Suppose now $M$ is a symplectic manifold and $\stackrel{M_0}{\stackrel{\downarrow}{\Phi}}$ is a Lagrangian fibration. Then Duistermaat's theory of action-angle
coordinates (see \cite{MG2}) implies that there is an action of the cotangent bundle $T^{\ast}\Phi$ on the fibers with a stabilizer lattice $\Lambda_b$. This of course induces a natural isomorphism $\xi: V_a/\Lambda_a \mapsto T^{\ast}\Phi / \Lambda_b$. There is also a dual isomorphism $\xi^{\ast}: T \Phi \mapsto V_a^{\ast}$ given explicitly by $\xi(v) = [i_v \omega]$, here $\omega$ is a symplectic structure and $v \in T\Phi$ is viewed as a normal vector field to an element of $\Phi$. Also the natural symplectic structure on $T^{\ast}\Phi$ projects to a symplectic structure on $T^{\ast}\Phi / \Lambda_b$ and hence on $V_a/\Lambda_a$.

If our fibration is a Special Lagrangian fibration then one can get a symplectic structure on the dual bundle $V_a^{\ast}/\Lambda_a^{\ast}$ as follows (this construction was done by Hitchin in \cite{Hit} and in coordinate-free way by Gross in \cite{MG2}):

We have a map $\alpha: V_a^{\ast} \mapsto T^{\ast}\Phi$ defined by periods of the closed form $Im \varphi$. Explicitly, let $u \in V_a^{\ast}$. We can view $u \in H^1(L,\mathbb{R})$. For $v \in T\Phi$ we define 
\begin{equation}\label{PD}
\alpha(u)(v) = [i_v Im\varphi] \cup u ([L])= [i_v Im\varphi](PD(u)) 
\end{equation}
Here $v$ is viewed as a normal vector field to $L$ and $PD(u)$ is a Poincare dual to $u$. One shows that for $u$ a section of $\Lambda_a^{\ast}$ (the integral cohomology lattice), $\alpha(u)$ is a closed 1-form on $\Phi$ and thus $\alpha$ induces a symplectic structure on $V_a^{\ast}/\Lambda_a^{\ast}$. This motivates the following definition 
\begin{dfn}
Let $\stackrel{M}{\stackrel{\downarrow }{\overline{\Phi}}}$ and 
$\stackrel{M'}{\stackrel{\downarrow }{\overline{\Phi}}}$ be 2 Special Lagrangian fibrations. Then $M'$ is a symplectic mirror to $M$ if the corresponding isomorphism $\rho: V_a^{\ast}/\Lambda_a^{\ast} \mapsto V_a'/\Lambda_a'$ over
$\Phi$ is a symplectomorphism.
\end{dfn}

To construct a topological mirror we make following observations: Let $F= 
\stackrel{W_a/ \Lambda_a}{\stackrel{\downarrow}{U}}$ be some torus fibration.
Let $a \in U$. Then we have a monodromy representation $\nu: \pi_1(U,a) 
\mapsto SL(W_a,\Lambda_a)$ (see \cite{MG1}). Moreover if $F' =
\stackrel{W_a' / \Lambda_a '}{\stackrel{\downarrow}{U}}$ is another fibration
and $K: W_a \mapsto W_a' $ is an intertwining isomorphism for monodromy 
representations, then $K$ induces
a natural fiberwise isomorphism between $F$ and $F'$. 

So we can try to compactify a dual fibration by trying to find local 
isomorphisms between it and the original fibration. Let $U$ be a neighbourhood
in $\Phi$ and $a \in U$.
Let $e_1, \ldots, e_n$ be
a basis for the lattice $\Lambda_a$ and $e^1, \ldots, e^n$ be a dual basis for the lattice $\Lambda_a^{\ast}$. Let $K: W_a/ \Lambda_a \mapsto W_a^{\ast}/ 
\Lambda_a^{\ast}$ be some 
linear map, which is given in terms of our bases by a matrix $K$. Let
$\alpha \in \pi_1(U,a)$ and $\nu(\alpha)$ be a monodromy map, which is given
in terms of a basis $(e_i)$ by a matrix $A$. Then a dual representation
$\xi^{\ast}(\alpha)$ on $W_a^{\ast}$ is given by a matrix $(A^T)^{-1}$ in a
basis $(e^i)$ (see \cite{MG1}).

So we need $K \cdot A = (A^T)^{-1} \cdot K$ i.e. $K= A^T \cdot K \cdot A$.

Now if $n=2$ there is a solution $K = \left( \begin{array}{clcr}
0 & 1 \\
-1 & 0 \\
\end{array} \right) $.

We return now to $M$. On a 6-torus $T^6$ we have a natural 
isomorphism between integral homology and cohomology of  SLag tori 
$T_{a,b,c}$. This isomorphism is invariant under $\mathbb{Z}_2
\oplus \mathbb{Z}_2$ action, hence it induces an isomorphism $\rho$ between
$\Lambda_a$ and $\Lambda_a^{\ast}$ outside of bad neighbourhoods.
Take a point $a$ on a boundary of a bad 
neighbourhood $Y$ in $\Phi$. Then because of the product structure of our 
fibration over $Y$ we see that the monodromy matrices of $V/ \Lambda$ in $Y$ 
look like $ \left( \begin{array}{clcr}
\ast & \ast & 0 \\
\ast & \ast & 0 \\
0 & 0 & 1 
\end{array} \right) $

It is clear from the above that the monodromy representation in Y is isomorphic to a dual representation.

So we can construct a topological mirror  $M'$ as follows: Let $\overline{W}$ be some
'bad'  neighbourhood in $M$ as before. Let $z_j=x_j+i \cdot y_j$ be
coordinates on a 6-torus. Then a map $\mu:T^6 \mapsto T^6$, $\mu(x_1,y_1,
x_2, y_2 , x_3,y_3)= (x_1,-y_2,x_2,y_1,x_3,y_3)$ commutes with $\mathbb{Z}_2^2$ action. Also $\mu$ maps the boundary $\partial 
\overline{W}$ to
itself. We take $\overline{W}$ and glue it to $M - \overline{W}$ by
$\mu$  and doing so for each 'bad' neighbourhood we obtain $M'$.

Now $\mu$ preserves SLag tori on $\partial \overline{W}$, and thus $M'$ naturally
acquires a structure of a fiber bundle over $\overline{\Phi}$ with generic fiber being a torus. We claim that $M'$ is a topological mirror of $M$.

Indeed we noted that outside of 'bad' neighbourhoods
there is a natural isomorphism $\rho$ between bundles $V_a$ and $V_a^{\ast}$
as before, and of course the bundle $V_a'$ of $M'$ is isomorphic to the 
bundle $V_a$ outside of 'bad' neighbourhoods, so $\rho$ can be viewed as an isomorphism between $\Lambda_a^{\ast}$ and $\Lambda_a'$. We want to extend $\rho$
inside of $\overline{W}$. First we need to check which isomorphism $\rho$ induces on $\partial \overline{W}$ via the gluing map $\mu$.

Let $L= T_{a,b,c}$ be a SLag torus contained in $\partial \overline{W}$. Let 
$z_j=x_j + i \cdot y_j$ be coordinates on $T^6$. Then $dy_1, \ldots, dy_3$ is a
basis for $H^1(L,\mathbb{Z})=\Lambda_a^{\ast}$ and $\partial_{y_1},\ldots, \partial_{y_3}$ is a dual basis for $H_1(L,\mathbb{Z})=\Lambda_a'$. Then 
\[\rho: dy_1 \mapsto -\partial_{y_2} ~ , ~ dy_2 \mapsto \partial_{y_1} ~ , ~
dy_3 \mapsto \partial_{y_3} \]   
As we saw, $\rho$ is an intertwining operator between the
monodromy representations on $V_a^{\ast}$ and $V_a'$. Hence $\rho$ extends to an isomorphism inside $\overline{W}$, and hence $M'$ is a topological mirror of $M$.

So far we viewed $M'$ just as a differential manifold (and one can easily show that $M'$ is diffeomorphic to $M$). We will see that
is has additional interesting structures.

Let $\omega'$,$\omega''$ and $\eta = Re(\eta)+ i \cdot Im(\eta)$ as before
(see p.14). Then we easily see that near $\partial \overline{W}$ the gluing map $\mu$ is an isometry. Also
$\mu^{\ast}(\omega')= Im \eta, \mu^{\ast}(Im \eta)= -\omega',\mu^{\ast}(Re \eta)
= Re \eta $.
Now
$Im(\eta) + (dx_3 \wedge dy_3)$ is a symplectic form on $\overline{W}$. So we
see that we can glue it to our symplectic form $ \omega' + (dx_3 \wedge 
dy_3)$ outside of $\overline{W}$ to get a symplectic form $\omega^{\ast}$ on $M'$.
Moreover near $\partial \overline{W}$, $\mu$ intertwines between $I$ and $K$ - almost 
complex structures defined by $\omega'$ and $Im \eta$. Thus we can glue
$K$ inside $\overline{W}$ to $I$ outside of $\overline{W}$ to get an almost complex structure
$I'$ on $M'$ compatible with $\omega^{\ast}$. We can glue a form $(Re \eta + i Im \eta) \wedge idz_3 $ outside of $\overline{W}$ to a form $(Re \eta -i \omega'') \wedge idz_3$ inside $\overline{W}$ to get a 
trivialization $\varphi'$ of a canonical bundle of $I'$.

A submanifold $L \in \Phi$, then viewed as a submanifold of $M'$, is  
Calibrated by $Re \varphi '$, which can be described by alternative 
conditions 
$\omega^{\ast} |_L=0, Im \varphi '|_L= 0$. So we can view our moduli-space on $M'$
as Special Lagrangian submanifolds, except for the fact that $I'$ is not
an integrable a.c. structure and so $M'$ is only symplectic. If we were
able to establish the fibration structure on $M$ by SLag submanifolds of the Calabi-Yau metric instead of $\omega'$,
we would have obtained a Calabi-Yau structure on the mirror. 

Finally we prove that $ \rho: V_a^{\ast}/\Lambda_a^{\ast} \mapsto V_a'/\Lambda_a'$ is a symplectomorphism and thus $M'$ is a symplectic mirror to $M$ according to Definition 3.3.2. The symplectic structure on $V_a^{\ast}/\Lambda_a^{\ast}$ was obtained from the map $\alpha: V_a^{\ast} \mapsto T^{\ast}\Phi$. Also the symplectic structure on $V_a'/ \Lambda_a' $ was obtained from a map $\xi': V_a' \mapsto T^{\ast}\Phi$. We will prove that
\[\alpha = \xi' \circ \rho \] and then we are done. This is obviously true outside of bad neighbourhoods. Let now $L$ be in one of bad neighbourhoods, so $L$
has a form $T \times S^1$. Let $\beta^1,\beta^2$ be (an oriented) basis for $H^1(T,\mathbb{Z})$. Then $\beta^1,\beta^2, [dy_3]$ is a basis for $H^1(L,\mathbb{Z})$, which we can view as a basis for $\Lambda_a^{\ast}$. Let $\beta_1, \beta_2, [S^1]$ be a corresponding dual basis for $H_1(L,\mathbb{Z})$, which we can also view as a basis for $\Lambda_a$ and $\Lambda_a'$. Then $\rho(\beta^1) = -\beta_2, \rho(\beta^2)= \beta_1 ~ and ~ \rho ([dy^3])= [S^1]$.

Let $v^i= \xi' (\beta_i)$ and $\gamma= \xi'([S^1])$ in $T^{\ast}\Phi$ (because of the product structure $\gamma$ can be viewed as a 1-form $dx_3$ on a moduli-space inside our bad neighbourhood).
Let $v_1,v_2, \partial_{x_3}$ be a dual basis of $T\Phi$. So if we view $v_i$
 as normal vector fields to the fiber then $i_{v_i}\omega^{\ast}$ represent cohomology classes $\beta^i$. But in a bad neighbourhood
 $\omega^{\ast} = Im\eta + dx_3 \wedge dy_3$, so
$[i_{v_i}Im\eta]= \beta^i$. 

Now by equation \ref{PD} we have \[ \alpha(\beta^i) (v_j)= [i_{v_j}Im\varphi] \cup \beta^i([L])= [i_{v_j}Im\eta \wedge -dy_3] \cup \beta^i (L)=\beta^j \cup -[dy_3]\cup \beta^i([L])= -\beta^i\cup \beta^j([T]) \] 
Also one easily shows that $\alpha(\beta^i)(\partial_{x_3})=0 $. So one deduces that
$\alpha(\beta^1)=-v^2 ~ ,~ \alpha(\beta^2)= v^1$ and also $\alpha([dy_3])= \gamma$. So $\alpha= \xi' \circ \rho $ and we are done. 

Remark: It is clear that applying these ideas we can get analogous results for
a Calabi-Yau 4-fold $N$ obtained from resolution of a quotient of an 8-torus by
$\mathbb{Z}_2^3$, there the generators of $\mathbb{Z}_2$ actions are

$\alpha: z_1, \ldots , z_4 \mapsto -z_1, -z_2, z_3, z_4 $

$\beta : z_1, \ldots , z_4 \mapsto  z_1, z_2, -z_3, -z_4$ 

and $\gamma: z_1, \ldots , z_4 \mapsto z_1, -z_2 + 1/2, -z_3 + 1/2 , z_4$

Indeed the resolution of the quotient by $\alpha$ and $\beta$ is a product of
2 $K^3$ surfaces with a product structure, there each $K^3$ has a metric 
Euclidean outside of bad neighbourhoods as before. 

For a fixed point set of $\gamma$ we introduce a bad neighbourhood $X$ in 
$z_2,z_3$ coordinates and consider a neighbourhood $Z = T^2 \times X \times
T^2$ in $T^8$. Then $\alpha$ and $\beta$ act freely on that neighbourhood.
Inside $Z$ we introduces structures as before and this way we get a Kahler
metric and a SLag torus fibration on $N$. 

\subsection{ Holomorphic functions near SLag Submanifolds }
In this section we examine holomorphic functions in a neighbourhood of a
special
Lagrangian submanifold. Such examples can be obtained for instance
from Calabi-Yau manifold $X$ in \(\mathbb{C}P^{n}\) defined as a zero locus 
of real
polynomials. In that case we have \( L = X \bigcap \mathbb{R}P^{n} \) a Special
Lagrangian submanifold. Let $P$ be some real polynomial of degree
$k$ without real roots. Then for any polynomial $Q$ of degree $k$ the
function \( \frac {Q}{P} \) is a holomorphic function on $X$ in a
neighbourhood of $L$. More generally let $L$ be a fixed point set of an 
antiholomorphic involution $\sigma$ and $h$ a meromorphic function on $M$.
Then obviously \( \overline{h \circ \sigma} \) is also a meromorphic function
on $M$ and so is \(g = h \cdot (\overline{ h \circ \sigma})+ 1 \). Also on $L$
$g$ is real valued and $\geq 1$. So $ f= 1/g$ is a holomorpic function in a 
neighbourhood of $L$. 

An immediate consequence from the fact that SLag submanifolds are `Special',
i.e. $ Im \varphi |_{L}= 0$ if the following 
\begin{thm}
: Let $L_0$ be Slag Submanifold and $f$ be a holomorphic
function
in a neighbourhood of $L_0$. Let $\xi$ be a function on our moduli-space,
\( \xi(L) = \int_{L} f \). Then $\xi$ is a constant function.
\end{thm}
{\bf Proof}: Consider the following $(n,0)$ form \( \mu = f \varphi\). Then
$\mu$ is holomorphic, hence closed and obviously $\xi(L)= \int_{L}\mu$. 
Q.E.D.\\

\noindent 
This yields a following corollary :\\
For \( 0 < \theta < \pi \) we denote by \( A_{\theta} \) an open
cone  in complex plane given by \( ( z=re^{i \rho } | r>0 , 0< \rho <
\theta ) \).
\begin{cor}
:  Let $M$ be a Calabi-Yau $n$-fold and $f$ a holomorphic
function on some domain $U$ in $M$. Let $L(t)$ be flow of  Slag
submanifolds contained in $U$ and \( p \in U \) a point s.t. \( f(p)=0\).
Suppose that the distance \(d(p, L(t)) \rightarrow 0 \) as \(t
\rightarrow \infty \). Then $L(t)$
cannot be contained in in the domain \( f^{-1} (A_{\theta}) \) for 
\( \theta < \frac{\pi}{2n} \). 
\end{cor}
Remark: This Corollary gives a restriction of how singular SLag currents
might look like .\\
{\bf Proof of Corollary 3.4.1}: Suppose $L(t)$ are contained in \(W= 
f^{-1}(A_{\theta})\) as above. We can find an \( \epsilon >0 \) s.t. 
\(g=f^{n+ \epsilon }\) is well defined an holomorphic on $W$ and 
\( g(W) \subset
A_{\frac{\pi}{2}} \). Then \( h= \frac{\pi}{2g} \) is holomorphic on $W$, 
\( h(W) \subset A_{\frac{\pi}{2}} \) and for $z$ close to $p$ we have \(|h(z)| \geq const \cdot 
d(z,p)^{-n- \epsilon} \).

Since \( \int_{L(t)} h \) is constant and \( Re(h),Im(h) >0\) on $L(t)$
then \( \int_{L(t)} |h| \) is bounded by a constant. \\
Take now any \( \delta >0 \) and pick $t$ and \( p_{t} \in L(t) \) s.t.
\(d(p,p_{t}) \leq \delta \). Consider \(B=B(p_{t}, \delta ) \bigcap L(t)
\). By Theorem 2.0.1, \( vol(B) \geq const \cdot \delta^{n} \) and on 
$B$ we have \( |h| \geq \frac{const}{\delta^{n+ \epsilon}} \). So \( \int_
{L(t)} |h| \geq \int_{B} |h| \geq const \cdot \delta^{- \epsilon}\).\\
Now $\delta$ was arbitrary - a  contradiction.      Q.E.D.\\

\noindent
Applying these ideas we can also get restriction on SLag submanifolds in
$\mathbb{C}^n$ which are asymptotic to a cone. We have the following theorem:
\begin{thm}
: Let \( L \subset \mathbb{C}^{n} \) be a special
Lagrangian submanifold
asymptotic to a cone $\Lambda$ and let \( z_{1} \ldots z_{n} \) be
coordinates on
$\mathbb{C}^n$. Then $L$ cannot be contained in the cone 
\[ B_{\theta}^{\delta}= ( (z_1, \ldots, z_n)| z_1 \in A_\theta ~,~ |z_1|>
\delta \cdot |z_i| ) \]
for $\delta > 0 ~,~ \theta < \pi/2n $.
\end{thm}
Remark : The order, to which $L$ is required to be asymptotic to a cone
will become clear from the proof.\\

\noindent 
{\bf Proof of the theorem}: Consider a flow $L(t)$ of SLag submanifolds in the unit ball
in
$\mathbb{C}^n$ with boundary in a unit sphere, \( L(t) = (z|t \cdot z \in L , |z| \leq  1 ) \). We wish to prove that \( \int_{L(t)} |z_{1}|^{-n-
\epsilon }  \) is uniformly bounded in $t$ for some $\epsilon > 0$ as in the proof of Corollary 3.4.1. This will lead us to a
contradiction as before because there are points in $L(t)$ which converge
to the origin in $\mathbb{C}^n$ for $t \rightarrow \infty$.

Let $d$ be the distance function to the origin on $L$. Let \(v =
\nabla d, w=\frac{v}{||v||^{2} } \). Then since $L$ is asymptotic to a
cone, $w$ will be a well-defined v. field outside some ball $B$ in $L$
and it's length will converge to 1 at $\infty$. We will also assume that
vector $w(x)$ is close to the line through $x$ and the origin, i.e. that 
there is a
function \( g : R_{+} \mapsto R_{+} \) s.t. the length of the orthogonal
component of $w(x)$ to this line is $\leq$ \( g(||x||) \) and \( \int_{[1,
\infty ) } \frac{g(t)}{t} dt < \infty \).

We extend $w$ inside
of $B$ to be a \( C^{\infty} \) v.field on $L$. Let \( \eta_{t} \) will
be flow of $w$ in time $t$, then the derivative of $d$ along \( \eta_{t}
\) is 1 outside $B$.\\
We can consider the corresponding flow \( \sigma_t \) on $L_s$,
\[ \sigma_{t}(x) = \frac{\eta_{t}(sx)}{s+t} ~,~ \sigma_{t} : L_{s}
\mapsto L_{s+t} \]
Let $v_s$ be a vector field on $L_s$ inducing the flow. One can
easily show that on the boundary of $L_s$ we have \( ||v_{s}|| \leq const
\cdot \frac{g(s)}{s} \). 

Pick $\epsilon$ small enough so that for $f(z) = z_{1}^{-n- \epsilon}$ we have $Ref,Imf $ are positive on $L_t$.
Let $h(t)= \int_{L_t}f=\int_{L_t}f \varphi$. We need to prove that $h(t)$ 
is a bounded function of $t$ and then we are done.  \\
Let $Q_t$ be the boundary of $L_t$. Then \[h'(t)=\int_{L_t}{\cal
L}_{v_t}f \varphi=\int_{L_t}d(i_{v_t}f \varphi)=\int_{Q_t}f \cdot i_{v_t}
\varphi \]
Conditions on $B_{\theta}^{\delta}$ imply that $|f|$ is uniformly 
bounded on $Q_t$. Also we know that \(|v_{t}| \leq \frac{g(t)}{t}\).\\
So \( |h'(t)| \leq const \cdot vol(Q_{t}) \frac{g(t)}{t} \). Now $Q_t$
converges to the base of the cone, so it's volume is bounded, so \(
|h'(t)| \) is an integrable function of $t$ by our assumptions, so    
$h(t)$ is uniformly bounded it $t$ .     Q.E.D. \\

\noindent 
Using those ideas for the 2-torus we get a following fact for analytic
functions of 1 complex variable :
\begin{lem}
: There is no holomorphic function $f$ from the open half disk \\
\(D = (re^{i \theta } | 0<r<1 , 0< \theta < \pi ) \) to itself s.t. 
\( |f(z)| \leq |z|^{2} \) .
\end{lem}
{\bf Proof}: Consider a domain \(D'=  (re^{i \theta} | 0<r<1 , 0< \theta <
\frac
{\pi}{2} ) \). \\
The map \( z \mapsto z^{2} \) is conformal from $D'$ to
$D$ and it is easy to see that is is enough to prove the claim for $D'$.
Also by iterating $f$ we can assume that \( |f(z)| \leq |z|^{4} \).

On a 2-torus \( T^{2} =\mathbb{R}^{2} (mod \mathbb{Z}^{2}) \) we consider SLag
submanifolds \[L_{t} = ((x,y)| y=t) \] We also have a Weierstrass
$\wp$-function, which has a pole of order 2 at the origin.

Consider
\(L_{\frac{1}{4}} \). We choose a constant $c$ so that \( P' =\wp +c \) satisfies \(
ReP' ~,~ImP' > 0 \) on \( L_{\frac{1}{4}} \). We look at  the flow 
\( L_{t} : t\rightarrow 0 \). We have 2 cases :\\
1) We have a (first) value $t_0$ s.t. \( ReP' (x, t_{0}) = 0 \) or
\(ImP'(x,t_{0}) = 0 \). W.l.o.g. we assume that the second case holds. Let
\(ReP'(x,t_{0}) = a \). Consider \( h = \frac{\pi}{2f(\sqrt{P'-
a})}  \).\\
Then as we approach $t_0$, h remains a holomorphic function with \( Reh ,
Imh > 0 \) and \( |h(z)| > const \cdot |z-(x+it_{0})|^{-2}  \). So as in 
Corollary 3.4.1 we get a contradiction.\\
2) If \(ReP' , Im P' \) remain positive as \( t \rightarrow 0 \) then we
get a contradiction by looking at \( \int_{L_{t}}|P'| \) .    Q.E.D. 

\section{Coassociative Geometry on a $G_2$ manifold}
For an oriented 7-manifold $M$, let $\bigwedge ^3 T^{\ast}M$ be a bundle of
3-forms on it. This bundle has an open sub-bundle $\bigwedge_+ ^3 M$ s.t.
$\varphi \in \bigwedge^3 T^{\ast}_p M$ is in $\bigwedge_+ ^3 M$ if there is a
linear isomorphism $\sigma: T_p M \mapsto \mathbb{R}^7$ s.t. $\sigma^{\ast}
\varphi _0 = \varphi$, there $\varphi_0$ is a standard $G_2$ 3-form on
$\mathbb{R} ^7$ (see \cite{Joy}, p. 294).

A global section $\varphi$ of $\bigwedge_+ ^3 M$ defines a topological $G_2$
structure on $M$. In particular this defines a Riemannian metric on $M$.
We will be interested in a closed $\varphi$. If $\varphi$ is also co-closed then
it is parallel and defines a metric with holonomy contained in $G_2$ (see
\cite{Joy}). In this case the form $\ast \varphi$ is a calibration and a
calibrated 4-submanifold $L$ is called a coassociative submanifold. This can
also be given by an alternative condition $\varphi|_L = 0$. For a closed
$\varphi$ we can define coassociative submanifolds by this condition. They no
longer will be Calibrated (because $\ast \varphi$ is not closed).
But nevertheless they are quite interesting because they admit an unobstructed
deformation theory. In fact we can copy a proof of theorem 4.5 in \cite{Mc} to show that their moduli-space $\Phi$ is smooth of dimension $b_2 ^+ (L)$
(proof of
that theorem in fact never used the fact that $\ast \varphi$ is closed).
If $L$ is coassociative as before and $p \in L$ we can identify the normal
bundle to $L$ at $p$ with self-dual 2-forms on $L$ by a map \[v \mapsto
i_v \varphi \] for $v$ a normal vector to $L$. Thus a tangent space to $\Phi$
at $L$ can be identified with closed self-dual 2-forms on $L$.

In a similar way to SLag Geometry we have a following lemma for finite group
actions on $M$ :
\begin{lem}
Suppose a finite group $G$ acts on $M$ preserving $\varphi$. Suppose $L$ is
a coassociative submanifold, $G$ leaves $L$ invariant and acts trivially on the second
cohomology of $L$. Then $G$ leaves every element of the moduli-space $\Phi$
through $L$ invariant.\
Moreover, suppose $g \in G$ and $x \in M - L$ is an isolated fixed point of
$g$. Then $x$ is not contained in any element in $\Phi$.
\end{lem}
The proof of this lemma is completely analogous to proof of Lemma 3.1.1.

Next we want to point out an example in which a $G_2$ manifold is a fibration
with generic fiber being a coassociative 4-torus. Our manifold $M$ is obtained from
resolution of a quotient of a 7-torus by a finite group. We hope to give a
systematic way of producing such examples in a future paper.

Let the group be $\mathbb{Z}_2^3$ with generators \[ \alpha: (x_1, \ldots, x_7)
\mapsto (-x_1,-x_2,-x_3,-x_4,x_5,x_6,x_7) \]
\[\beta : (x_1,\ldots,x_7) \mapsto (-x_1+1/2,1/2-x_2,x_3,x_4,-x_5,-x_6,x_7) \]
\[\gamma : (x_1,\ldots,x_7) \mapsto (-x_1, x_2, 1/2-x_3, x_4, -x_5,x_6,-x_7)\]

(compare with \cite{Joy}, p.302). We will follow Joyce's exposition of that
example.

The fixed point locus of each generator is a disjoint union of 3-tori. Their
fixed loci don't intersect and their compositions have no fixed points. Around
each fixed point the quotient looks like $V=T^3 \times B/\pm1$, there $B$ is a
ball in $\mathbb{R}^4$. We will show how to get a $G_2$ structure on resolution of singularities $\overline{V}$. We will treat a fixed locus of, say, $\alpha$. 

Let $x_1, \ldots, x_4$ be coordinates on $\mathbb{R}^4$. Let
$\omega_1,\omega_2,\omega_3$ be a standard HyperKahler package on
$\mathbb{R}^4$. 
For coordinates $x_5,x_6,x_7$, let $\delta_i$ be a dual 1-form to $x_{8-i}$. 
Then the $G_2$ 3-form on $\mathbb{R}^7$ is \[ \varphi= \omega_1 \wedge \delta_1 + \omega_2 \wedge \delta_2 + \omega_3 \wedge \delta_3 + \delta_1 \wedge \delta_2 \wedge \delta_3 \]
We can use either one of 3 complex structures on $\mathbb{R}^4$ to identify it with $\mathbb{C}^2$. That way for each singularity of the form $\mathbb{R}^4/\pm1$ we get a resolution of a singularity $\overline{U}$ and a (non-integrable) HyperKahler package on $\overline{U}$ in a similar way to section 3.3.
Suppose, for example, that we used a complex structure $I$ on $\mathbb{R}^4$. Then $\omega_2$ and $\omega_3$ lift to $\overline{U}$. Also $\omega_1$ is replaced by a
Kahler form $\omega_1'$ (as in Section 3.3).
 On $\overline{V}= \overline{U} \times T^3$
we can consider \[ \varphi'= \omega_1' \wedge \delta_1 + \omega_2 \wedge \delta_2 +\omega_3 \wedge \delta_3 + \delta_1 \wedge \delta_2 \wedge \delta_3 \]
and $\varphi'$ will be a closed $G_2$ form. Doing so for each fixed locus we get a $G_2$ structure on $M$. D. Joyce proved that this 3-form can be deformed to a parallel $G_2$ form. We will use the form $\varphi'$ because we can construct a coassociative 4-torus fibration on $M$ for $\varphi'$.

On $T^7$ we have a following coassociative 4-torus fibration \[T_{a,b,c}=
((x_1,\ldots,x_7)| x_1=a, x_3=b, x_6=c) \]
Note that the 4 coordinates on each $T_{a,b,c}$ are chosen so that each
generator of $\mathbb{Z}_2^3$ acts non-trivially on exactly 2 of those coordinates. Those $T_{a,b,c}$ become coassociative tori on $M$ and fill a big 
open neighbourhood of $M$. We would like to see what happens then
these tori enter a 'bad' neighbourhood $V$ above. For that we have to consider a bigger tubular neighbourhood $W=X \times T^3$, there
\[ X= (v \in T^4| v = w + (0,a,0,b), ~there~ w \in U ~and~ a,b \in \mathbb{R}/ \mathbb{Z} ) \]
So $W$ are exactly those points in $T^7$ which are contained on
$T_{a,b,c}$ for a torus $T_{a,b,c}$, which intersects $V$. We have a resolution of singularities $\overline{W}$, which can be viewed as a neighbourhood in $M$.

We would like to investigate coassociative tori in $\overline{W}$. As we mentioned, we have 3 different ways to resolve a singularity using either one of the
structures $I,J,K$. We have 2 different cases :

1) The structure we are using is either $I$ or $K$.  We can assume that it is $I$. The $G_2$ form looks like \[\omega_1' \wedge \delta_1 + \omega_2 \wedge \delta_2 + \omega_3 \wedge \delta_3 + \delta_1 \wedge \delta_2 \wedge \delta_3 \]
Our tori will look like 
\[T_c \times L \]
there $T_c$ is a torus in $x_5,x_6,x_7$ coordinates defined by condition $x_6=c$ and $L$ is in $\overline{U}$ defined by conditions: $\omega_1'|_L=0, \omega_3|_L=0 $, i.e it is a Special Lagrangian torus as in section 3.3. The results of
section 3.3 precisely apply to show that the compactified moduli-space fibers
the neighbourhood $\overline{W}$. 

2) The structure we are using is $J$. Then we have a package $\omega_1, \omega_2',\omega_3$ and the $G_2$ form looks correspondingly. Then our tori again look like $T_c \times L$, there $L$ satisfies $\omega_1|_L=0, \omega_3|_L = 0$, i.e.
they are holomorphic tori with respect to a structure $J$. So what we get is
a neighbourhood in $K^3$ with a standard holomorphic fibration over a neighbourhood in $S^2$. So in any case our manifold $M$ fibers with generic fiber being a coassociative 4-torus.

Finally we want to construct a topological mirror for this torus fibration. We will use definitions from section 3.3. It is clear that because of the local product structure the local monodromy representation is isomorphic to the dual representation.
Hence we can construct a dual fibration by performing a surgery for each 'bad'
neighbourhood $\overline{W}$ in a similar way to section 3.3. So for instance if  $\overline{W}$ is a neighbourhood of $\alpha$ then we glue $\overline{W}$ to $M - \overline{W}$ along a boundary by a map \[ \eta: (x_1,\ldots,x_7) \mapsto (x_1,-x_4,x_3,x_2,x_5,x_6,x_7) \]
Also near the boundary of $\overline{W}$ we have $\eta^{\ast}(\omega_1)= \omega_3 ~,~ \eta^{\ast}(\omega_3) = -\omega_1 ~,~ \eta^{\ast}(\omega_2)= \omega_2 $.
So we can glue a closed $G_2$ form $\varphi^{\ast}= \omega_3 \wedge \delta_1 + \omega_2 \wedge \delta_2 - \omega_1 \wedge \delta_3 + \delta_1 \wedge \delta_2 \wedge \delta_3$ inside of $\overline{W}$ to $\varphi'$ outside of
$\overline{W}$ to get a closed $G_2$ form $\varphi''$ on the mirror. Also the mirror fibration is a fibration by coassociative 4-tori with respect to $\varphi''$.  

Remark: The original example in Joyce's paper \cite{Joy} was a quotient by a slightly different $\mathbb{Z}_2^3$ action with generators \[ \alpha: x_1, \cdots , x_7 \mapsto -x_1,-x_2,-x_3,-x_4,x_5,x_6,x_7 \] \[ \beta : x_1 , \cdots , x_7 \mapsto 
-x_1,1/2-x_2,x_3,x_4,-x_5,-x_6,x_7 \] \[ \gamma: x_1, \cdots , x_7 \mapsto 1/2 - x_1,x_2,1/2-x_3,x_4,-x_5,x_6,-x_7 \]
We get a closed $G_2$ form $\varphi$ on the resolution of singularities similarly to previous example. Our manifold $M$ again will be fibered by coassociative 4-tori if we start from a family $T_{a,b,c} = ((x_1, \ldots, x_7)|x_1=a,x_2=b,x_7=c)$. Indeed we consider a neighbourhood $U_i$ of a fixed component of one of
the generators and a bigger neighbourhood $X_i=(v \in T^7|v= u + (0,0,a_1,a_2,a_3,a_4,0) ~ s.t.~ u \in U ~ and ~ a_i \in \mathbb{R}/\mathbb{Z}) $. Then $X_i$ are disjoint, so we can use  product structure on $X_i$ to get a coassociative torus fibration on $M$ as in the previous example.

Massachusetts Institute of Technology 

E-mail : egold@math.mit.edu

\end{document}